\documentclass [11pt]{article}

\usepackage[dvips,pagebackref,colorlinks=true,linkcolor=blue,citecolor=blue,pdfstartview=FitH]{hyperref}

\usepackage{amssymb,amsmath}

\newcommand\Z{\mbox{$\mathbb Z$}}

\newcommand\R{\mbox{$\mathbb R$}}
\newcommand\C{\mbox{$\mathbb C$}}
\newcommand\F{\mbox{$\mathbb{F}$}}

\def\B{\{0,1\}}
\def\H{\{-1,1\}}
\def\J{[-1,1]}

\def\np{{\rm NP}}
\def\pcp{{\rm PCP}}

\newlength{\pgmtab}  
 \newenvironment{program}{%
\begin{tabbing}\hspace{0em}\=\hspace{0em}\=%
\hspace{\pgmtab}\=\hspace{\pgmtab}\=\hspace{\pgmtab}\=\hspace{\pgmtab}\=%
\hspace{\pgmtab}\=\hspace{\pgmtab}\=\hspace{\pgmtab}\=\hspace{\pgmtab}\=%
\+\+\kill}{\end{tabbing}}

\newcommand {\ACCEPT}{{\bf accept}}

\def\poly{{\rm poly}}

\let\eps=\varepsilon
\def\epsilon{\varepsilon}
\def\phi{\varphi}
\def\bx{{\bf x}}
\def\bg{{\bf g}}
\def\bzero{{\bf 0}}

\newcommand\ip[2]{\langle #1 \rangle_{U^{#2}}}
\newcommand\lip[2]{\langle #1 \rangle_{LU^{#2}}}

\newtheorem{theorem}{Theorem}
\newtheorem{lemma}[theorem]{Lemma}

\newtheorem{proposition}[theorem]{Proposition}

\newtheorem{remark}{Remark}
\newtheorem{definition}[theorem]{Definition}
\newenvironment{proof}{\noindent {\sc Proof:}}{$\Box$ \medskip}

\newcommand\E{\mathop\mathbb{E}}
\newcommand\pr{\mathop{\bf {Pr}}}
\newcommand\var{\mathop{\bf Var}}
\newcommand\I{{\rm  Inf}}
\newcommand\XI{{\rm  XInf}}
\newcommand{\ignore}[1]{}

\def\({\left(}
\def\){\right)}

\def\lemmaref#1{\hyperref[#1]{Lemma~\ref*{#1}}}

\title{Gowers Uniformity, Influence of Variables, and PCPs}
\author{Alex Samorodnitsky\thanks{Institute of Computer Science, Hebrew University.
{\tt salex@huji.ac.il}} \and Luca Trevisan\thanks{Computer Science Division, U.C. Berkeley.
{\tt luca@cs.berkeley.edu}. This paper
is based upon work supported by the National Science Foundation under grant CCF 0515231
and by the US-Israel Binational Science Foundation under grant 2002246.}}

\begin {document}

\maketitle

\sloppy

\begin{abstract}
Gowers \cite{Gowers98,Gowers01}
introduced, for $d\geq 1$, the notion of {\em dimension-$d$ uniformity} $U^d(f)$ of a function $f: G\to \C$, where
$G$ is a finite abelian group.
Roughly speaking, if a function has small Gowers uniformity of dimension $d$, then it ``looks random'' 
on certain structured subsets of the inputs.

We prove the following ``inverse theorem.'' Write $G=G_1 \times \cdots \times G_n$
as a product of groups. 
If a bounded balanced function $f:G_1 \times \cdots G_n \to \C$ is such that $U^{d} (f)\geq \epsilon$,
then one of the coordinates of $f$ has {\em influence} at
least $\epsilon / 2^{O(d)}$. Other inverse theorems are known \cite{Gowers98,Gowers01,GT:inverse,S05}, and
$U^3$ is especially well understood, but the properties of
functions $f$ with large $U^d(f$), $d\geq 4$, are not yet well
characterized.
 
The  {\em dimension-$d$ Gowers inner product} $\ip{ \{ f_S \} }{d}$
of a collection $\{ f_S \}_{S\subseteq [d]}$ of functions is a related measure of pseudorandomness. 
The definition is such that if all the functions $f_S$
are equal to the same fixed function $f$, then $\ip{ \{ f_S \} }{d} = U^d(f)$. 

We prove that if $f_S:G_1\times \cdots \times G_n \to \C$ 
is a collection of bounded functions such that $|\ip{ \{ f_S \} }{d} |\geq \epsilon$ and
at least one of the $f_S$ is balanced,
then there is a variable that has influence at least $\epsilon^2 / 2^{O(d)}$ for at least
four functions in the collection. 

Finally, we relate the acceptance probability of the ``hypergraph long-code test'' proposed
by  Samorodnitsky and Trevisan  to the Gowers inner product of the functions being tested
and we deduce the following result: if the Unique Games Conjecture is true, then for every $q\geq 3$
there is a PCP characterization of NP where the verifier makes $q$ queries, has almost perfect
completeness, and soundness at most $2q/2^q$. For infinitely many $q$, the soundness is $(q+1)/2^q$.
Two applications of this results are that, assuming that the unique games conjecture is true,
it is hard to approximate Max $k$CSP within a factor $2k/2^k$ (or even $(k+1)/2^k$, for infinitely many $k$),
and it is hard to approximate
Independent Set in graphs of degree $D$ within a factor $(\log D)^{O(1)} / D$. 
\end{abstract}

\section{Introduction}

We return to the study of the relation between number of queries and error probability
in probabilistically checkable proofs.

The PCP Theorem \cite{AS92,ALMSS92} states that it is possible
to encode certificates of satisfiability for SAT instances (and, more
generally, for every problem in NP) in such a way that a probabilistic verifier
can check the validity of the certificate with high confidence after
inspecting only a constant number of bits. We write $SAT \in \pcp_{c,s}[r(n),q]$
if there is a verifier that uses at most $r(n)$ random bits, where $n$ is the size of the formula,
accepts encoding of valid proofs with probability at least $c$ (the {\em completeness} probability
of the verifier) and accepts purported encodings of proofs with probability at most $s$
(the {\em soundness} error of the verifier) if the formula is unsatisfiable. 
The PCP Theorem states that there exists a constant $k$ such that $SAT \in \pcp_{1,1/2}[O(\log n),k]$. 
Improvements and variants of the PCP Theorem and their applications to the study
of the approximability of optimization problems are too many to summarize here,
and we refer the reader to the chapter on hardness of approximation in
Vazirani's book \cite{V01} and to some recent survey papers \cite{A02:icm,F02:icm,T04:approx}.

In this paper we are interested in the following question: for a given number of queries,
what is the highest confidence that we can have in the validity of the proof? That is, 
for a given value of $q$, what is the smallest value $s=s(q)$ for which $SAT\in \pcp_{1-\delta,s+\delta}[O(\log n),q]$
for every $\delta>0$?  We call this parameter $s$ the {\em soundness} of the PCP construction.
A good measure of the trade-off between the number $q$ of queries and the soundness $s$
is the so-called {\em amortized query complexity}, defined as $\bar q = q/ (\log_2 s^{-1})$.

A simple argument shows that, unless $P=NP$, $s$ must be at least $1/2^q$, that is,
the amortized query complexity must be at least $1$. A more
careful argument gives a lower bound of $2/2^q$ \cite{T96} on the soundness, which was
recently improved to $\Omega(\frac q {\log q} \cdot \frac1 {2^q})$ by Hast \cite{H05}.
(Hast's result can also be stated as giving a lower bound of $1+ (1-o(1))\frac {\log q}{q}$ to the amortized
query complexity of a $q$-query PCP.)
The PCP Theorem shows that we can have $s= 1/2^{\Omega(q)}$, and
 the authors showed that can have $s\leq 2^{2\sqrt q}/2^{q}$ \cite{ST00}.
 (That is, the amortized query complexity can be as low as $1+ O(1/\sqrt q)$.)
Our proof was simplified by H\aa stad and Wigderson \cite{HW01}, and the
soundness was  improved to $s\leq 2^{\sqrt {2q}}/2^{q}$ by  Engebretsen and Holmerin \cite{EH05}.
As we discuss below, $2^{\Theta(\sqrt q)}/2^q$ was a natural limit for the soundness
achievable with current techniques.

In this paper, assuming Khot's Unique Games Conjecture \cite{K02:unique}, we
present an improvement to $s= (q+1)/2^q$. Our analysis is based on a theorem, which
is probably of independent interest, bounding the Gowers uniformity of a given
function in terms of the influence of its variables.

\subsection{Linearity Tests and PCP}

The {\em linearity testing problem} captures most of the technical difficulties
of the the construction of PCP constructions, and it is a good starting point.

Let us call a function $f:\B^n \to \H$ {\em linear} if it
is a homomorphism between the group $\B^n$ (where the operation is bitwise XOR,
written additively) and
the group $\H$ (where the operation is multiplication). Equivalently,
$f$ is linear if and only if it can be written as $f(x_1,\ldots,x_n) = (-1)^{\sum_{i\in S} x_i}$
for some set $S\subseteq [n]$. We use the notation $\chi_S(x) :=  (-1)^{\sum_{i\in S} x_i}$.

In the linearity testing problem we are given oracle access to 
a boolean function $f: \B^n \to \H$ and we would like to distinguish
between the following extreme settings:

\begin{enumerate}
\item $f$ is linear;
\item for every $S$, the agreement between $f$ and $\chi_S$ is  at most $1/2 + \epsilon$.
\end{enumerate}

By {\em agreement} between a function $f$ and a function $g$ we mean the fraction of inputs
on which they are equal. We say that a test has {\em error-probability} at most $e$ if
in case (1) it accepts with probability 1 and in case (2) it accepts with probability
at most $e+ \epsilon'$, where $\epsilon'\to 0$ when $\epsilon \to 0$.

Blum, Luby and Rubinfeld \cite{BLR90} define a very simple such test, that makes only
three queries into $f$:

\noindent\fbox{\begin{minipage}{2in}
\begin{program}
BLR-Test (f)\+\\
 choose $x,y$ uniformly at random in $\B^n$\\
 \ACCEPT\ if and only if $f(x)\cdot f(y) = f(x+y)$ 
\end{program}\end{minipage}}

\bigskip

Bellare et al. \cite{BCHKS95} give a tight analysis of this test, showing
that if it accepts with probability at least $1/2+\epsilon$, then 
$f$ has agreement at least $1/2+\epsilon$ with some linear function.
According to our above definition, the BLR test has error probability at most 1/2.

There are at least two ways in which such a result needs to be extended before
it can be used in a PCP construction. 

First of all, we would like to consider
a case where two or more functions are given as an oracle, and the test wants
to distinguish between the following cases:

\begin{enumerate}
\item The functions given as an oracle are all equal to the same linear function
\item No two functions have agreement more than $1/2+\epsilon$ (or less than $1/2-\epsilon$)
with the same linear function
\end{enumerate}

There is a natural extension of the BLR test to this setting:

\noindent\fbox{\begin{minipage}{2in}
\begin{program}
3-functions-BLR-Test $(f,g,h)$\+\\
 choose $x,y$ uniformly at random in $\B^n$\\
 \ACCEPT\ if and only if $f(x)\cdot g(y) = h(x+y)$ 
\end{program}\end{minipage}}

Aumann et al. \cite{AHRS01} show that if this test accepts with probability $1/2+\epsilon$,
then there is a linear function $\chi_S$ such that $f,g,h$ have all agreement at least $1/2 + \epsilon/3$
with either $\chi_S$ or $-\chi_S$.

The second change is that, for the sake of PCP constructions, we are especially interested
in linear functions $\chi_S$ with a small $S$. We call functions of the form $\chi_{ \{ i \} } (x) = (-1)^{x_i}$
{\em long codes}. In a {\em long code test} we are given several functions and, in a commonly used
definition, we want to distinguish the following cases:

\begin{enumerate}
\item The functions given as an oracle are all equal to the same long code;
\item For every small $S$, no two functions have agreement more than $1/2+\epsilon$ (or less than $1/2-\epsilon$)
with  $\chi_S$. 
\end{enumerate}

We say that such a test has error probability at most $e$ if, whenever it accepts
with probability more than $e+\epsilon$, there are constants $\epsilon',d$, depending only on $\epsilon$, and
a set $S$, $|S|\leq d$, such that at least two of the given functions have
agreement at least $1/2+\epsilon'$ with $\chi_S$ or $-\chi_S$.

A test satisfying this definition can be obtained from the BLR test by adding noise
to each query. Let $\mu_\delta$ be the probability distribution over $\B^n$
defined by picking $n$ independent copies of a biased coin that returns 0 with probability
$1-\delta$ and 1 with probability $\delta$, and consider the following test:

\noindent\fbox{\begin{minipage}{2in}
\begin{program}
$\delta$-noisy-3-functions-BLR-Test $(f,g,h)$\+\\
 choose $x,y$ uniformly at random in $\B^n$\\
 sample $\eta_1,\eta_2,\eta_3$ indipendently according to $\mu_\delta$\\
 \ACCEPT\ if and only if $f(x+\eta_1)\cdot g(y+\eta_2) = h(x+y+\eta_3)$ 
\end{program}\end{minipage}}

\bigskip

It is easy to see that in case (1), that is, when $f,g,h$ are equal to the same
long code, then the test accepts with probability at least $(1-\delta)^3$.
H\aa stad \cite{H97} shows\footnote{Such an analysis is implicit in
H\aa stad's paper, and the result, as stated here, appears explicitely in \cite{AHRS01}.} 
that if the test accepts with probability at least $1/2+\epsilon$,
then there is a set $S$ of size at most $\poly(\epsilon^{-1},\delta^{-1})$ such that
$f,g,h$ all have agreement at least $1/2 + \epsilon/3$ with $\chi_S$ or $-\chi_S$.

Having a $q$-query test for this problem with error probability at most $e$
is enough to construct a PCP characterization of NP with query complexity $q$
and soundness error about $e$, provided that Khot's {\em unique games conjecture}
\cite{K02:unique} is true. (Further refinements are needed to derive an unconditional
result in \cite{H97}.)

Given this logical path from the basic linearity testing problem to the task of constructing
PCPs, our plan in \cite{T98:stoc,ST98,ST00} was to devise  linearity tests with
good trade-offs between number of queries and error probability, and then
``translate'' such  tests into a PCP construction.

In \cite{ST00} we devise a linearity tester whose asymptotic trade-off between number of queries
and error probability is optimal. The test, for $k\geq 2$ is defined as follows:

\noindent\fbox{\begin{minipage}{2in}
\begin{program}
Complete-Graph-Test $(f,k)$\+\\
 choose $x_1,\ldots,x_k$ uniformly at random in $\B^n$\\
 \ACCEPT\ if and only if\+\\ for every $i\neq j$, $f(x_i)\cdot f(x_j) = f(x_i+x_j)$ 
\end{program}\end{minipage}}

Note that the test has query complexity $q=k+ {k\choose 2}$ and runs ${k \choose 2}$
correlated copies of the BLR test. If $f$ is linear, then it is clearly accepted
with probability 1. If $f$ has agreement at most $1/2 +\epsilon$, then we already know that
each of the ${k \choose 2}$ tests accepts with probability at most $1/2+\epsilon$.
In \cite{ST00} we show that the ${k\choose 2}$ tests behave almost mutually independently, and the
probability that all accept is at most $1/2^{k\choose 2} + \epsilon' \approx 1/2^{q-\sqrt{2q}}$.
We also extended the test to the noisy case, the case of several functions, and the
setting (which we  do not describe in this paper) which is sufficient to derive
an unconditional PCP characterization of NP.\footnote{In this last step, we lost something in the
soundness error, which became
 $1/2^{q-2\sqrt{q}}$, where $q$ is the number of queries. This was recently improved 
 to $1/2^{q-\sqrt{2q}}$, the same bound of the basic linearity test, by Engebretsen and Holmerin \cite{EH05}.}
 
One might have thought that the following test would have achieved an even better trade-off
between number of queries and error probability:

\noindent\fbox{\begin{minipage}{2in}
\begin{program}
Complete-Hypergraph-Test $(f,k)$\+\\
 choose $x_1,\ldots,x_k$ uniformly at random in $\B^n$\\
 \ACCEPT\ if and only if\+\\ for every $S\subseteq [k]$: $|S|\geq 2$. $\prod_{j\in S} f(x_j)= f(\sum_{j\in S} x_j)$ 
\end{program}\end{minipage}}

In the hypergraph test we make $q=2^k-1$ queries and run $2^{k}-k-1$ correlated copies
of the BLR test. If the tests behaved almost mutually independently for $f$ far from linear,
then the test would have error probability $\approx (q+1)/2^q$.

Unfortunately, the above statement fails already when $k=3$ if we consider
the function $f(x) = (-1)^{x_1x_2 + \cdots + x_{n-1}x_n}$. Such a function is
very far from linear, but the 7 conditions of the hypergraph test for $k=3$ are not almost
independent. In fact, much more generally, we prove in \cite{ST00} that if we run a $q$-query hypergraph test
on $f$, then the test accepts with probability at least $2^{-q+\Omega(\sqrt q)}$. 
In Section \ref{sec:lowerbound}
we generalize this result and prove that any linearity test that makes $q$ queries
and that accepts linear functions with
probability $\geq c$ must accept $f$ with probability at least $(1-c) + 2^{-q+\Omega(\sqrt q)}$.

Even though there is no more room for improvement in the basic linearity testing
problem (or its multi-function version, which is only harder), there is still hope
for PCP, by using a more relaxed formulation of the long code test.
Several recent papers \cite{DS02,K02:unique,KR03,KKMO04,CKKRS05,KV05} define
a long code test based not on closeness to linear functions $\chi_S$ with small
$S$ but based on the notion of 
degree-$d$ influence. We will define such a notion later (Section \ref{sec:prelim}).
For now, it suffices to say that if we can  solve a certain (relaxed) variant
of the linearity test with a given query complexity and error probability, then we 
can also test the long code with the same query complexity and error probability.
We say that the relaxed test has error probability $e$ if:

\begin{enumerate}
\item if $f$ is linear, then the test accepts with probability 1;
\item if the test accepts with probability $e+\epsilon$, then there is a coordinate $i$
that has influence at least $\delta(\epsilon)$ for $f$.
\end{enumerate}

Influence (as opposed to ``degree-$d$ influence'') has a simple definition
for boolean functions: the influence of $i$ for $f$ is the probability
that $f$ is a non-constant function after we randomly fix all coordinates
of $f$ except the $i$-th. If $f$ has agreement $1/2+\epsilon$ with a non-constant
linear function $\chi_S$, then $f$ has variables of influence at least $2\epsilon$,
but there are functions $f$ that have influential variables even though they
are very far from all linear functions (this is why this test is a relaxation).
Intriguingly, in the function $f(x)= (-1)^{x_1x_2+\cdots+x_{n-1}x_n}$ all
variables have influence $1/2\pm o(1)$.

As we shall see later, we prove that the complete hypergraph test
has error probability only $1/\left( 2^{2^k-k-1} \right)$ with respect
to the above definition. That is, if a function $f$ is
accepted by the test with probability more than $1/\left( 2^{2^k-k-1} \right) +\epsilon$,
then one of the variables of $f$ has influence at least $\delta(\epsilon)$. We
then generalize the test to the setting of several function (in which case
we show that at least two functions have an influential variable in common)
and to the case of testing the long code, resulting in a conditional PCP
construction assuming the unique games conjecture.

\subsection{{S}zem\'eredi's Theorem and Gowers Norms}

We will use a definition that has been introduced by Gowers \cite{Gowers98,Gowers01}
in his seminal papers on a new proof of Szem\'eredi's Theorem. 

Szem\'eredi's Theorem states that any subset of the integers of positive density contains
arbitrarily long arithmetic progressions. The quantitative statement is that for
every $\delta,k$ there is a $n(\delta,k)$ such that if $A$ is an arbitrary subset
of the integers $\{1,\ldots,N\}$, $|A|\geq \delta N$, and $N\geq n(\delta,k)$, then
$A$ must contain an arithmetic progression of length $k$, that is, there are integers $a,b$ such
that $a,a+b,\ldots,a+(k-1)b$ all belong to $A$.

In Szem\'eredi's original proof, which introduced the famous Szem\'eredi Regularity Lemma,
$n(\delta,k)$ grows {\em very} fast with $k$ and $\delta^{-1}$:
it is a tower of exponentials whose height is a tower of exponentials whose height
is polynomial in $k$ and $\delta^{-1}$. 

The case $k=3$ had been settled earlier by Roth, with a simpler analytical
proof that gave a reasonable (doubly exponential) relation between $n(\delta,3)$ and $\delta^{-1}$.
Roth's proof (like all other proofs about arithmetic sequences in dense sets) is based on an 
iterative argument on $\delta$: if $\delta$ is a sufficiently large constant (say, 90\%), then
a random length-3 progression is contained in $A$ with positive probability, and so
$A$ contains some length-3 progressions. Otherwise (blurring, for the moment, the
distinction between progressions mod $N$ and true progressions), one writes the
fraction of length-3 progressions in $A$ as the number one would expect if $A$
where a random set of density $\delta$, that is, a $\delta^3$ fraction of all length-3 progressions,
minus an error term that equals $\sum_{g\neq 0} \hat f_A^3(g)$, the summation of the cubes
of the non-zero Fourier coefficients of $f_A$, the characteristic function of $A$.
One then considers two cases: if all Fourier coefficients of $f_A$ are small,
then the error term is smaller than $\delta^3$, and so $A$ contains a positive
fraction of all length-3 progressions in $\Z_N$. If $f_A$ has a large Fourier
coefficient, then one can reduce the task of finding a length-3 progression in $A$
to the task of finding a length-3 progression in a certain subset $A'$ of $\{1,\ldots,N'\}$
of density $\delta' > \delta + \Omega(\delta^2)$. In this reduction, $N'$ is about $\sqrt{N}$,
one does no more than $O(\delta^{-1})$ such reductions, so it's enough to 
start with $N= 2^{2^{O(1/\delta)}}$. We mention these technical details because
they are very similar to the analysis of the basic BLR linearity test in \cite{BCHKS95},
where the acceptance probability of the verifier is written as $\frac 12 + \frac 12 \sum_\alpha \hat f^3(\alpha)$,
and so (i) if all Fourier coefficients of $f$ are small, then $\sum_\alpha \hat f^3(\alpha)$ is
small, and the test accepts with probability close to $1/2$, while (ii) if one Fourier
coefficient is large, then $f$ is close to a linear function (this part is trivial).

A way to look at both proofs is to think of a function $f$ as being ``pseudorandom''
if all its Fourier coefficients are small, and of a set $A\subseteq [N]$ to be pseudorandom
if all the non-zero Fourier coefficients of its characteristic function are small. Then
one proves that, for a pseudorandom function, the values $f(x),f(y),f(x+y)$ are nearly
independent random bits, and so, in particular, $f(x)f(y)=f(x+y)$ happens with probability
approximately 1/2. For a psedorandom set of density $\delta$, the events $x\in A$, $(x+r)\in A$,
$(x+2r)\in A$ also behave nearly independently for random $x,r$, and they happen simultaneously
with probability approximately $\delta^3$.

To prove Szem\'eredi's Theorem for $k\geq 4$ one might try to show that a pseudorandom set,
as defined above, has approximately the expected number of length-$k$ progression. This,
unfortunately, does not seem to be true, and sets whose characteristic function is defined
in terms of a degree-2 polynomial are basic counterexamples even for $k=4$. (See \cite{Gowers01}.) Recall
that something similar happens in the hypergraph test, where a function defined in terms
of a degree-2 polynomial is very far from linear (and so all its Fourier coefficients are small),
but the tests performed in the hypergraph test do not behave independently.

Gowers \cite{Gowers98,Gowers01} resolves this problem by introducing a more refined
measure of pseudorandomness. For $d\geq 1$ and $f:G\to \R$, where $G$ is a group, he
defines the quantity
\[ U^d(f) := \E_{x,x_1,\ldots,x_d} \left[ \prod_{S\subseteq [d]} f\left( x + \sum_{i\in S} x_i \right) \right] \]
and there is a similar definition for $f:G\to \C$, in which all evaluations of $f()$ when $|S|$ is odd
are conjugated. (See definitions in Sections \ref{sec:prelim} and \ref{sec:complex}.)

Gowers goes on to prove that if $A\subseteq \Z_N$ is a subset of density $\delta$,
$f:Z_N\to [-1,1]$ is a normalized version of the characteristic function, and $U^d(f)$
is sufficiently small (as a function of $d$ and $\delta$, but not of $N$) then 
$A$ must contain arithmetic progressions of length $d+1$. The {\em hard} part
of Gowers's work is then to show that if $U^d(f)$ is large, then finding progressions
in $A$ reduces to finding progressions in a set $A'\subseteq \Z_{N'}$ of density
strictly larger than $\delta$. 

Towards this goal, Gowers proves certain structural properties of functions
$f:\Z_N \to \C$, $N$ prime, having non-trivially large $U^d$ value. Such
functions are shown to have a certain ``local correlation'' with degree $(d-1)$-polynomials.
Green and Tao \cite{GT:inverse} study functions $f:G\to \C$ with large $U^3$, and, provided that the
order of $G$ is not divided by and 2 and 3, prove a certain ``global correlation''
between such functions and degree-2 polynomials.\footnote{The ``globality'' of this result deteriorates
with the order of the group, and a result as stated is false for, say, $\Z_N$, $N$ prime.} 
Samorodnitsky \cite{S05} proves such a result for functions $f:\B^n \to \R$.

Not much is known about functions $f:G\to \C$ having large $U^d$ when $d\geq 4$
and $G$ is a general group.

\subsection{Our Results}

We prove that if $f:G_1\times \cdots \times G_n \to \C$ is a bounded balanced
function and $U^d(f)\geq \epsilon$, then there is a variable that has
influence at least $\epsilon / 2^{O(d)}$ for $f$. Above, we only defined
influence in the boolean case, but a more general definition applies
to functions mapping from an arbitrary product space into the complex numbers.

Green and Tao introduce a notion related to the $U^d$ measures of pseudorandomness.
For $2^d$ functions $\{ f_S \}_{S\subseteq [d]}$, $f_S:G\to \R$, their 
{\em Gowers inner product} is defined as
\[ \ip{ \{ f_S \} }{d} := \E_{x,x_1,\ldots,x_d} \left[ \prod_{S\subseteq [d]} f_S\left( x + \sum_{i\in S} x_i \right) \right] 
\]
In the case of complex-valued functions, the terms with odd $|S|$ are conjugated. Note
that if all the functions $f_S$ are identical to the same function $f$ then
$\ip{ \{ f_S \} }{d} = U^d(f)$.

Our second result is that if $f:G_1\times \cdots \times G_n \to \C$ are bounded functions
such that at least one of them is balanced, and
$\ip{ \{ f_S \} } d \geq \epsilon$, then there is a variable that has
influence at least $\epsilon^2 / 2^{O(d)}$ for at least {\em four} of
the functions in the collection.

Finally, we come back to the analysis of the hypergraph test.
H\aa stad and Wigderson \cite{HW01} significantly simplified
the analysis of the graph test of \cite{ST00} by using an
averaging argument  that reduces the analysis of the graph
test to the analysis of the 3-functions BLR test.

We apply a similar averaging argument and reduce the analysis
of the hypergraph test to the task of bounding expressions
of the form

\begin{equation} \label{eq:lineargowers}
\E_{x_1,\ldots,x_d} \prod_{S\subseteq [d]} f_S\left( \sum_i x_i \right) \end{equation}

where $d\leq k$ and the functions $f_S:\B^n \to \H$ are derived in a certain
way from the function being tested. The expression looks very similar to the
expression for the Gowers inner product, and in fact it is easy to see that
if the Expression in (\ref{eq:lineargowers}) is large, then the Gowers inner
product of a related set of functions is also large. From our results, it follows
that at least four of these new functions must share an influential variable,
from which it follows that two of $f_S$ must share an influential variable and
so the function being tested has an influential variable. (See Section \ref{sec:lowerbound}.)

This analysis easily extends to the case of testing multiple functions and to the
noisy case. We  present the analysis of the noisy, multi-function case
in Section \ref{sec:pcp}.

This leads,  under the unique games conjecture, to a PCP characterization of NP
with $q$ queries and error probability $(q+1)/2^q$, if $q$ is of the form $2^k-1$.
As corollary, we obtain $(q+1)/2^q$ hardness of approximation for $qCSP$ and
$(\poly\log D)/D$ hardness of approximation for independent set in graphs of
maximum degree $D$. In terms of amortized query complexity, we show
that $1+ (1+o(1))\frac {\log q}{q}$ is achievable assuming the unique games conjecture,
and Hast's algorithm~\cite{H05} implies that this is tight up to the lower order term.

\subsection{Organization of the paper} We develop the basic theory of
influence of variables and Gowers uniformity for the case $G= \Z_2$ in
Section \ref{sec:prelim}. We prove our connection between influence
of variables and Gowers uniformity in the case $G=\Z_2^n$ in
Sections \ref{sec:prod-inf}, \ref{sec:uniformity}, \ref{sec:ip}.
Section \ref{sec:complex} shows how to generalize our results
to the case of general abelian groups; only the result of Section \ref{sec:prod-inf}
requires a new proof, while the other results follow by making syntactic changes in the
proof for the boolean case.   An application to generalized linearity testing
is presented in Section \ref{sec:lowerbound}, together with lower bounds 
establishing the tightness of our analysis.
An application to PCP is
presented in Section \ref{sec:pcp}.

The paper ends up being quite long because, for the sake of readability, we
first prove some results in important special cases before proceeding to
the more general statements. A reader interested only in the PCP application
may  skip Section \ref{sec:complex}. 
A reader interested only in our results on Gowers uniformity may skip
Sections \ref{sec:lowerbound} and \ref{sec:pcp}.

\section{Preliminaries}
\label{sec:prelim}

In this section we develop the theory of Gowers uniformity and Fourier transforms
for functions $f:\Z_2^n \to \R$. Generalized definitions for
the setting of functions $f:G_1\times \cdots \times G_n \to \C$ will be given later in Section \ref{sec:complex}.

For a positive integer $n$, we use $[n]$ to denote the set $\{1,2,\ldots,n\}$. We use the convention
that $[0]$ is the empty set.

For 
two reals $a<b$ we use $[a,b]$ to denote the interval $\{ x\in \R : a\leq x\leq b\}$.

\subsection{Fourier Analysis}

For a subset $S\subseteq [n]$, define the function $\chi_S:\B^n \to \R$ as
\[ \chi_S(x_1,\ldots,x_n) = (-1)^{\sum_{i\in S} x_i} \]
We use the convention that an empty sum equals zero, so $\chi_\emptyset(x)=1$ for every $x$.

For two functions $f,g :\B^n \to \R$, define their  inner product as
\[ \langle f, g\rangle := \E_{x\in \B^n} f(x)g(x) \]

Then it is easy to see that the functions $\chi_S$ are {\em orthonormal} with respect to this
inner product, that is, for every $S$,
\[ \langle \chi_S , \chi_S \rangle = 1\]
and, for every $S\neq T$,
\[ \langle \chi_S , \chi_T\rangle = 0\]
This implies that the functions $\chi_S$ are linearly independent and, since there 
are $2^n$ such functions and the space of functions $f:\B^n \to \R$ has dimension $2^n$,
it follows that the $\chi_S$ are a {\em basis} for such space.

Every function $f:\B^n \to \R$ can therefore be written as
\[ f(x) = \sum_{S\subseteq [n]} \hat f(S) \chi_S(x) \]
where the coefficients $\hat f(S)$ in the linear combination satisfy
\[ \hat f(S) = \langle f, \chi_S\rangle \]
and are called the {\em Fourier coefficients} of $f$. The function $\hat f()$ mapping
sets $S$ into coefficients $\hat f(S)$ is the {\em Fourier transform} of $f$. We will
make use of the following equality, called {\em Parseval's identity} (or {\em Plancherel's identity}).

\begin{equation} \sum_S \hat f^2 (S) = \E_x f^2 (x) \end{equation}

In particular, if $f:\B^n \to [-1,1]$, then $\sum_S \hat f^2 (S) \leq 1$, and so $|\hat f(S)| \leq 1$
for every $S$.

\subsection{Influence of Variables}

If $f:\B^n \to \H$ is a boolean function, then the influence of $i$ for $f$ is defined as
\begin{equation} \label{eq:def:boolean-influence}
 \I_i(f) := \pr_x [ f(x) \neq f(x+e_i)] 
\end{equation}
where $e_i\in \B^n$ is the vector that has a 1 in the $i$-th position and zeroes everywhere else.
That is,  $\I_i(f)$ is the probability that, if we pick a random $x$, the value of $f$
at $x$ is different from the value of $f$ at the point that we get by flipping the $i$-th bit of $f$.

It is easy that see that $\I_i(f)$ satisfies the following identities.

\begin{equation} \label{eq:alternate-influence}
 \I_i(f) = \frac 14 \E_x ( f(x) - f(x+e_i) )^2 =\sum_{S : i\in S} \hat f^2(S)
\end{equation}

For a general real-valued function $f:\B^n \to \R$, we will define the influence of $i$ for $f$ as
\[ \I_i(f):= \frac 14 \E_x ( f(x) - f(x+e_i) )^2 =\sum_{S : i\in S} \hat f^2(S) \]
Note that Equation \ref{eq:def:boolean-influence} is not valid if $f$ is a general real-valued function.

We will make use of the following simple fact, that is valid for every function
$f:\B^n \to \R$: if $S \neq \emptyset$ and $i\in S$, then
\[ \I_i(f) \geq \hat f^2 (S) \]
and, in particular, for every $f:\B^n \to \R$,
\begin{equation} \label{eq:influence:simple}
 \max_i \I_i (f) \geq \max_{S \neq \emptyset} \hat f^2(S) 
\end{equation}

For boolean functions $f:\B^n \to \H$, the stronger inequality $\max_i \I_i (f) \geq \max_{S\neq \emptyset} |\hat f(S)|$
holds, but we will not use it.

For a function $f:\B^n \to \R$, a threshold $d\leq n$, and an index $i\in [n]$, we
define the {\em degree-$d$ influence of $i$} for $f$ as
\[ \I^{\leq d}_i (f):= \sum_{S: i\in S, |S|\leq d} \hat f^2 (S) \ . \]

We will make use of the following simple fact: if $f$ is a bounded function, then
not too many variables can have large low-degree influence. Specifically:

\[ \sum_{i=1}^n \I^{\leq d}_i(f) = \sum_{i=1}^n \sum_{S: i\in S, |S|\leq d} \hat f^2 (S) = 
\sum_{S: i\in S, |S|\leq d} |S| \hat f^2(S) \leq d \sum_S \hat f^2(S) \]
In particular, if $f:\B^n \to \J$, then $\sum_{i=1}^n \I^{\leq t}_i(f) \leq d$,
and so at most $d/\epsilon$ variables can have degree-$d$ influence larger than $\epsilon$.

\subsection{Cross-Influence}

For the application to PCP, the following definition will be useful.
Let ${\cal F} = f_1,\ldots,f_k$ be a collection of (not necessarily distinct) functions $f_j : \B^n \to \R$. Then 
the {\em cross-influence} of $i$ for $\cal F$ is defined as
\begin{equation} \XI_i (f_1\ldots,f_k) := \max_{j\neq h} \min \{ \I_i(f_j) , \I_i(f_h) \} \end{equation}

In other words, $\XI_i (f_1,\ldots,f_k) \geq \epsilon$ if and only if there are at least two 
functions $f_j,f_h$, with $j\neq h$, such that $\I_i(f_j)\geq \epsilon$ and $\I_i(f_h) \geq \epsilon$.
Conversely, $\XI_i (f_1,\ldots,f_k) \leq \epsilon$ if there is at most one function
$f_j$ such that $\I_i (f_j) > \epsilon$.

More generally, for a threshold  $t$, we define the $t$-cross influence of $i$ for $\cal F$ as
\begin{equation}
t\XI_i(f_1,\ldots,f_k) := \max_{ j_1,\ldots,j_t \in [k], {\rm\ all\ distinct} }\min
\{ \I_i(f_{j_1}) , \ldots, \I_i(f_{j_t}) \}
\end{equation}
That is, $t\XI_i ({\cal F})$ is 
the smallest $\epsilon$ such that there are at least $t$ functions $f_{j_1},\ldots,f_{j_t}$ in $\cal F$ 
such that coordinate $i$ has influence at least $\epsilon$ for all of them. Note that
$\XI({\cal F}) = 2\XI({\cal F})$.

If all the functions $f_j$ are equal to the same function $f$, then 
$k\XI_i(f_1,\ldots,f_k) = \XI_i(f_1,\ldots,f_k) = \I_i (f)$. 

Similarly, we define the degree-$d$ cross-influence of $\cal F$ as
\begin{equation} 
\XI^{\leq d} _i (f_1\ldots,f_k) := \max_{j\neq h} \min \{ \I^{\leq d}_i(f_j) , \I^{\leq d}_i(f_h) \} \end{equation}

\subsection{Gowers Uniformity}

\begin{definition}[Gowers Uniformity]
Let $f:\B^n \to \R$ be a function, and $d\geq 1$ be an integer. The {\em dimension-$d$ uniformity}
of $f$ is defined as

\[ U^d(f) := \E_{x,x_1,\ldots,x_d} \prod_{S\subseteq [d]} f\left(x+\sum_{i\in S} x_i\right) \]
\end{definition}

\begin{remark} Here we use a terminology and notation that is a hybrid between the one
of Gowers \cite{Gowers98,Gowers01} and the one of Green and Tao \cite{GT:primes,GT:inverse}. 
What we call {\em dimension-$d$ uniformity} is called {\em degree-$(d-1)$ uniformity}
by Gowers, and no notation is introduced for it. Gowers also introduces the notation
$|| f ||_d$, which equals, in our notation, $\( U^d (f) \)^{1/2^d}$. Gowers proves that
$|| \cdot ||_d$ is a norm, and he does not give it a name. 
Green and Tao use the notation $|| f||_{U^d}$ for $\( U^d (f) \)^{1/2^d}$,
and call it the {\em (dimension-$d$) Gowers norm}.
\end{remark}

Here are expressions for the first few values of $d$:

\begin{eqnarray*}
 U^1(f) & = & \E_{x,y} f(x)f(x+y) = \left( \E_x f(x) \right)^2 \\
 U^2(f) & = & \E_{x,y,z} f(x)f(x+y)f(x+z)f(x+y+z) \\
& = & \E_y \left (\E_x f(x)f(x+y) \right)^2 \\
 U^3(f) & = & \E_{x,y,z,w} f(x)f(x+y)f(x+z)f(x+y+z)f(x+w)f(x+y+w)f(x+z+w)f(x+y+z+w)\\
&  = & \E_{y,z} \left ( \E_x  f(x)f(x+y)f(x+z)f(x+y+z) \right)^2 
\end{eqnarray*}

The above examples suggest the use of the following notation. 

For a function
$f:\B^n \to \R$ and elements $x_1,\ldots, x_d\in \B^n$, define

\[ f_{x_1,\ldots,x_d} (x) := \prod_{S\subseteq [d]} f\left(x+ \sum_{i\in S} x_i \right) \]
Then we have
\[ U^d(f) = \E_{x,x_1,\ldots,x_d} f_{x_1,\ldots,x_{d}} (x) = \E_{x_1,\ldots,x_{d-1}} 
\left( \E_x f_{x_1,\ldots,x_{d-1}} (x)
\right)^2 \]

Define the {\em dimension-$d$ Gowers inner product} of a collection
$\{ f_S \}_{S\subseteq [d]}$  of (not necessarily different) functions $f_S:\B^n \to \R$,
as

\begin{equation}
\ip{| \{ f_S \} }{d} := \E_{x,x_1,\ldots,x_{d}} \left[ \prod_{S\subseteq [d]} f_S \left( x + \sum_{i\in S} x_i \right) \right] \end{equation}

Note that, in particular, if all the functions $f_S$ are equal to the same function $f$ then
$\ip{\{ f_S \}}{d} = U^d (f)$.

\subsection{Unique Games}

A unique game \cite{K02:unique} is a constraint satisfaction problem such that every constraint
is of the form $y=f_{x,y} (x)$, where $x,y$ are variables ranging over a finite set $\Sigma$, which 
we call the {\em alphabet},
specified as part of the input, and $f_{x,y} : \Sigma \to \Sigma$  is a permutation. 
Given a unique game, we are interested in finding the assignment of values to the variables
that satisfies the largest number of constraints.

More formally, a unique game is a tuple $(V,E,\Sigma,\{ f_{x,y} \}_{(x,y)\in E})$
where $V$ is a set of variables, $E$ is a set of pairs of variables (corresponding to constraints),
and, for every $(x,y)\in E$, the function $f_{x,y}:\Sigma \to \Sigma$ is a permutation.
Note that $(V,E)$ is a graph, which we call the {\em constraint graph} of the unique game.
We want to find an assignment $A: V\to \Sigma$ that maximizes the number of satisfied constraints,
that is, the number of pairs $(x,y)\in E$
such that $A(y) = f_{x,y} (A(x))$. The {\em value} of an assignment is the fraction
of constraints satisfied by the assignment; the value of a unique game is the 
value of an optimum assignment.

For example, the following is a unique game with $V=\{v_1,v_2,v_3,v_4\}$
and $\Sigma=\{ a,b,c\}$:

\begin{eqnarray*}
v_3 & = & {a\ b\ c \choose c\ b\ a} (v_1)\\
v_3 & = & {a\ b\ c \choose a\ c\ b} (v_2)\\
v_1 & = & v_2\\
v_4 & = & {a\ b\ c \choose b\ c\ a} (v_2)
\end{eqnarray*}

Where we use the notation $\left( \begin{array}{ccc} a & b & c\\ f(a) & f(b) & f(c) \end{array} \right)$
to represent a function $f: \{ a,b,c\} \to \{a,b,c\}$.
The reader can verify that the value of the above unique game is $3/4$.

The {\em unique games conjecture} is that for every $\gamma > 0$ there is a $\sigma=\sigma(\gamma)$ such that it is NP-hard
to distinguish unique games of value $\geq 1-\gamma$ from unique games of
value $\leq \gamma$, even when restricted to instances where $|\Sigma| \leq \sigma$ and where
the constraint graph is bipartite.

For our application, we will need a variant of unique games, that we call $d$-ary unique game.
In a $d$-ary unique game, a constraint is specified by a $d$-tuple $v_1\ldots,v_d$ of variables
and a $d$-tuple of permutations $f_1,\ldots,f_d :\Sigma\to \Sigma$. An assignment $A:V\to \Sigma$
{\em strongly satisfies} the constraint if $f_1(A(v_1)),\ldots,f_d(A(v_d))$ are
all equal; an assignment {\em weakly satifies} the constraint it if $f_1(A(v_1)),\ldots,f_d(A(v_d))$ are
not all different. 

The following result is a rephrasing of a result by Khot and Regev \cite{KR03}.

\begin{theorem} \label{th:outer}
If the unique games conjecture is true, then for every $d$ and every $\gamma$
there is a $\sigma=\sigma(d,\gamma)$ such that, given a $d$-ary unique game with alphabet
size $\sigma$, it is NP-hard to distinguish
the case in which there is an assignment that strongly satisfies at least  a $1-\gamma$
fraction of constraints from the case where every assignment weakly satisfies
at most a $\gamma$ fraction of constraints.
\end{theorem}

We define the {\em strong value} of an assignment to a $d$-ary unique game as 
the fraction of constraints that are strongly satisfied by the assignment.
The strong value of a unique game is the largest strong value among all assignments.
The {\em weak value} of an assignment and of  a unique game are similarly defined.
Note that the weak value is always at least as large as the strong  value.

\section{Influence of Product of Functions}
\label{sec:prod-inf}

In this section we prove a bound on the influence of a function
of the form $f(x):= f_1(x) \cdot f_2(x) \cdots f_k(x)$ in terms
of the influence of the functions $f_j$. Such a bound will be
very useful in the proofs of our main results.

In the boolean case, the bound is just a simple union bound.

\begin{lemma}\label{lm:influence-shift}
Let $f_1,\ldots,f_k:\B^n \to \H$ be boolean functions, and 
define $f(x) = f_1(x) \cdot f_2(x) \cdots f_k (x)$.

Then, for every $i\in [n]$
\begin{equation} \label{eq:influence-shift}
 \I_i(f) \leq \sum_j  \I_i (f_j)
\end{equation}
\end{lemma}

\begin{proof}
Using the formula for the influence of boolean functions, we see that
\begin{eqnarray*}
 \I_i(f) & = & \pr_x [ f(x) \neq f(x+e_i)] \\
& = &  \pr_x [ f_1(x)\cdots f_k(x)
\neq f_1(x+ e_i )\cdots f_k(x +e_i) ]  \\
& \leq & \sum_{j=1}^k \pr_x [ f_j(x) \neq f_j(x+ e_i) ]\\
& = & \sum_j \I_i(f_j)
\end{eqnarray*}
\end{proof}

For general real-valued functions, we cannot hope to achieve
the nice bound of Equation \ref{eq:influence-shift}. 
Suppose, for example, that $n=1$ and that all functions $f_j$
are defined as follows: $f_j(0) = 1-\epsilon$, $f_j(1) = 1$.

Then, we have $I_1(f_j) = \frac 14 \epsilon^2$. When we
define $f(x) := \prod_j f_j(x)$, we get $f(0)=(1-\epsilon)^k$
and $f(1)=1$, and so $I_i(f) = \frac 14 (1-(1-\epsilon)^k)^2$,
which is about $\frac 14 k^2 \epsilon^2$ for small $\epsilon$, or
about $k$ times the sum of the influences of the functions $f_j$.
The following Lemma achieves such a tight bound.

\begin{lemma}\label{lm:influence-shift-general}
Let $f_i:\B^n \to \J$ be $k$ functions, and define
\[ f(x) := f_1(x)\cdots f_k(x) \]
Then, for every $i\in [n]$, $\I_i(f) \leq k \cdot \sum_j \I_i (f_j)$
\end{lemma}

\begin{proof}
We begin by proving the following claim:
\begin{equation}\label{eq:cute-inequality}
 \forall a_1,\ldots,a_k,b_1,\ldots,b_k \in \J. ~~~\left| \prod_i a_i - \prod_i b_i \right| \leq  \sum_i \left|a_i - b_i \right| 
\end{equation}

We prove Inequality \ref{eq:cute-inequality} by first expressing the right-hand side
as a telescoping sum of a sequence of ``hybrids,'' and then by using the 
triangle inequality and the fact that all $a_i$
and $b_i$ have absolute value at most 1.

\begin{eqnarray*}
&& | a_1\cdot a_2 \cdots a_k  - b_1 \cdot b_2 \cdots b_k|\\
 & = & |(a_1\cdot a_2 \cdots a_k) - (b_1\cdot a_2 \cdots a_k)
+(b_1\cdot a_2 \cdots a_k) - (b_1\cdot b_2 \cdot a_3\cdots a_k) +(b_1\cdot b_2 \cdot a_3\cdots a_k) -\\ && \cdots
+ (b_1\cdots b_{k-1}\cdot a_k) -(b_1\cdots b_k) |\\
&=& | (a_2\cdots a_k)\cdot (a_1-b_1) + (b_1\cdot a_3 \cdots a_k) \cdot (a_2 - b_2) 
+ \cdots + (b_1\cdots b_{k-1}) \cdot (a_k-b_k) |\\
& \leq & |a_2\cdots a_k|\cdot |a_1 - b_1| + |b_1\cdot a_3 \cdots a_k| \cdot |a_2 - b_2| 
+ \cdots + |b_1\cdots b_{k-1}|\cdot |a_k-b_k |\\
& \leq & |a_1 - b_1| + |a_2 - b_2|  + \cdots + |a_k-b_k |
\end{eqnarray*}

If we square both sides and apply Cauchy-Schwartz, we get
\begin{eqnarray*}
\left ( \prod_i a_i - \prod_i b_i \right)^2 & \leq & \left ( \sum_i |a_i - b_i | \right)^2\\
& \leq &  k \cdot \sum_i (a_i - b_i)^2
\end{eqnarray*}
 
To summarize our progress so far, we have proved the following claim:
\begin{equation}\label{eq:ab}
 \forall a_1,\ldots,a_k,b_1,\ldots,b_k \in \J. ~~~\left( \prod_i a_i - \prod_i b_i \right)^2 \leq  k\cdot
\sum_i \left(a_i - b_i \right)^2 
\end{equation}

We are now ready to prove the Lemma. For every $x$, using Inequality \ref{eq:ab}, we have
\[ ( f(x) - f(x+e_i))^2 =  \left(\prod_i f_i (x) - \prod_i f_i(x+e_i) \right)^2
\leq  k \cdot \sum_j (f_j(x) - f_j(x+e_i) )^2 \]
and of course the same inequality remains valid if we take the average over $x$, so we have
\[ \I_i(f) = \frac 14 \E_x (f(x) - f(x+e_i))^2 \leq k \frac 14 \E_x \sum_j (f_j(x) - f_j(x+e_i))^2
\leq k \sum_j \I_i(f_j) \] 
\end{proof}

\section{Low Influence Implies Small Gowers Uniformity}
\label{sec:uniformity}

We are going to show that, for balanced functions $f:\B^n \to \H$, if $U^d(f)$ is large, then
one of the variables of $f$ has high influence.

\begin{lemma}\label{lm:boolean-influence-one}
Let $f:\B^n \to \H$ be a function and $d\geq 1$ be an integer. Then
\[  U^d(f) \leq U^1(f) + (2^{d-1}-1) \max_i \I_i (f) \]
\end{lemma}

\medskip

\begin{proof}
The case $d=1$ is trivial. Let $d\geq 2$, and define $I:= \max_i \I_i(f)$. We will prove
\begin{equation}
U^d(f) \leq U^{d-1}(f) + 2^{d-2} I
\end{equation}
which immediately implies the statement of the lemma.

We write
\[ U^d(f) = \E_{x_1,\ldots,x_{d-2}}\left [ \E_{x,y,z} f_{x_1,\ldots,x_{d-2}}(x)
f_{x_1,\ldots,x_{d-2}}(x+y)f_{x_1,\ldots,x_{d-2}}(x+z)f_{x_1,\ldots,x_{d-2}}(x+y+z) \right] \]
 \[ = \E_{x_1,\ldots,x_{d-2}}\sum_\alpha \hat f_{x_1,\ldots,x_{d-2}}^4(\alpha )\]
\[ = \E_{x_1,\ldots,x_{d-2}}\sum_{\alpha\neq \emptyset} \hat f_{x_1,\ldots,x_{d-2}}^4(\alpha )
+ \E_{x_1,\ldots,x_{d-2}} \hat f_{x_1,\ldots,x_{d-2}}^4(\emptyset ) \]

We separately bound the two terms in the last expression. 

For every $a_1,\ldots,a_{d-2} \in \B^n$, using 
Lemma \ref{lm:influence-shift}, we get that, for every $i\in \alpha$,
$\I_i(f_{a_1,\ldots,a_{d-2}}) \leq 2^{d-2} \I_i (f) \leq 2^{d-2} I$
and so we have that, for every $a_1,\ldots,a_{d-2} \in \B^n$,

\[ \sum_{\alpha\neq \emptyset} \hat f_{a_1,\ldots,a_{d-2}}^4(\alpha )
\leq \max_{\alpha\neq \emptyset} \hat f_{a_1,\ldots,a_{d-2}}^2(\alpha ) \leq  \max_i \I_i (f_{a_1,\ldots,a_{d-2}})
\leq 2^{d-2} I \]
and so
\[ \E_{x_1,\ldots,x_{d-2}}\sum_{\alpha\neq \emptyset} \hat f_{x_1,\ldots,x_{d-2}}^4(\alpha ) \leq 2^{d-2} I \]

Regarding the other term,
\[ \E_{x_1,\ldots,x_{d-2}} \hat f_{x_1,\ldots,x_{d-2}}^4(\emptyset ) \leq 
\E_{x_1,\ldots,x_{d-2}} \hat f_{x_1,\ldots,x_{d-2}}^2(\emptyset )
 = \E_{x_1,\ldots,x_{d-2}} \left( \E_x f_{x_1,\ldots,x_{d-2}}(x )\right)^2 = U^{d-1}(f) 
\]
\end{proof}

The same argument also applies to general bounded real-valued functions.

\begin{lemma}\label{lm:boolean-influence-two}
Let $f:\B^n \to \J$ be a function and $d\geq 1$ be an integer. Then
\[  U^d(f) \leq U^1(f) + 4^d \max_i \I_i (f) \]
\end{lemma}

\medskip

\begin{proof} 
Define $I:=\max_i \I_i (f)$. It suffices to prove that, for $d\geq 2$, 
\begin{equation}
U^d(f) \leq U^{d-1}(f) + 2^{2d-4} \max_i \I_i (f)
\end{equation}
We repeat the proof of \lemmaref{lm:boolean-influence-one} verbatim,
except that we use \lemmaref{lm:influence-shift-general}
instead of \lemmaref{lm:influence-shift} to get an upper
bound for $\I_i(f_{a_1,\ldots,a_{d-2}})$. Because of the worse
bound in \lemmaref{lm:influence-shift-general}, we only get the bound
\[   \I_i(f_{a_1,\ldots,a_{d-2}}) \leq 2^{2d-4} \I_i (f) \leq 2^{2d-4} I \]
and the rest of the proof proceeds with no change, except for the term $2^{2d-4}$
instead of $2^{d-2}$.
\end{proof}

We remark that our bound for the boolean case is nearly tight.

\begin{lemma}
For every fixed  $d\geq 2$ and large $n$, there is a function $f:\B^n \to \H$
such that $U^{1} (f) = o_n(1)$, $U^{d} (f)= 1$ and $\max_i \I_i (f) \leq \frac 1 {2^{d-2}} +o_n(1)$.
\end{lemma}

\begin{proof}
Consider the function $$ f(x_1,\ldots,x_n) := (-1)^{x_1 x_2 \cdots x_{d-1} + x_d \cdots x_{2d-2} + \cdots}$$
\end{proof}

\section{Low Cross-Influence Implies Small Gowers Inner Product}
\label{sec:ip}

The main result of this section is that if a collection of functions has small cross-influence,
then it has small Gowers inner product, provided that at least one of the functions is balanced.

\begin{lemma} \label{lm:main}
For $d\geq 2$, let $\{ f_S \}_{S\subseteq [d]}$ be a collection
of functions $f_S: \B^n \to [-1,1]$ such that
\begin{itemize}
\item $f_{[d]}$ is balanced, that is, $\E_x f_{[d]}(x) = 0$;
\item $4\XI_i (\{ f_S \}) \leq \epsilon$ for every $i$.
\end{itemize}

Then 
\[ \ip{\{ f_S \} }{d} \leq \sqrt{\epsilon} \cdot 2^{O(d)} \]
\end{lemma}

Before proving \lemmaref{lm:main}, we establish a variant of
a result of Aumann et al. \cite{AHRS01} which will be
useful in the inductive step of the proof of \lemmaref{lm:main}.

\begin{lemma}\label{lm:fourfunctions}
For every four bounded functions $f_1,f_2,f_3,f_4: \B^n \to [-1,1]$,
\[ \left | \sum_\alpha \hat f_1(\alpha) \hat f_2(\alpha) \hat f_3(\alpha) \hat f_4 (\alpha) \right|
\leq 4\max_{\alpha} \min \{ |\hat f_1(\alpha) |,|
 \hat f_2(\alpha)|, | \hat f_3 (\alpha)|,| \hat f_4 (\alpha)|\} \]
 \end{lemma}

\begin{proof}
Let \[\epsilon:= \max_\alpha  \min \{ |\hat f_1(\alpha) |,|
 \hat f_2(\alpha)|, | \hat f_3 (\alpha)|,| \hat f_4 (\alpha)|\}\]
For $i=1,\ldots,4$, let $S_i$ be the family of all sets $\alpha$ such that $|\hat f_i(\alpha)|\leq \epsilon$.
By definition, the union of the families $S_i$ contains all subsets $\alpha \subseteq [n]$. We can then write

\begin{eqnarray*}
 &  & \left| \sum_\alpha \hat f_1(\alpha) \hat f_2(\alpha) \hat f_3(\alpha) \hat f_4 (\alpha) \right| \\
& \leq & \sum_\alpha |\hat f_1(\alpha)| \cdot | \hat f_2(\alpha)| \cdot | \hat f_3(\alpha)| \cdot | \hat f_4 (\alpha)|\\
& \leq & \sum_{i=1}^4 \sum_{\alpha\in S_i} |\hat f_1(\alpha)| \cdot | \hat f_2(\alpha)| \cdot | \hat f_3(\alpha)| \cdot | \hat f_4 (\alpha)|\\
& < &  \epsilon \left(
\sum_{\alpha \in S_1} | \hat f_2(\alpha)| \cdot | \hat f_3(\alpha)| \cdot | \hat f_4 (\alpha)|
+\sum_{\alpha \in S_2} | \hat f_1(\alpha)| \cdot | \hat f_3(\alpha)| \cdot | \hat f_4 (\alpha)| \right. \\
&& \left. +\sum_{\alpha \in S_3} | \hat f_1(\alpha)| \cdot | \hat f_2(\alpha)| \cdot | \hat f_4 (\alpha)|
+\sum_{\alpha \in S_4} | \hat f_1(\alpha)| \cdot | \hat f_2(\alpha)| \cdot | \hat f_3 (\alpha)| \right)\\
& \leq & 4\epsilon
\end{eqnarray*}
The last inequality follows from the fact that for every three functions $f,g,h:\B^n \to [-1,1]$ we
have
\[ \sum_\alpha | \hat f(\alpha) | | \hat g(\alpha) | | \hat h(\alpha) |
\leq \sum_\alpha | \hat f(\alpha) | | \hat g(\alpha) | \leq \sqrt{\sum_\alpha \hat f^2(\alpha)}  
\sqrt{\sum_\alpha \hat g^2(\alpha)} \leq 1
\]
\end{proof}

We proceed with the proof of our main result of this section.

\bigskip

\begin{proof}[Of \lemmaref{lm:main}]
We want to prove 
\begin{equation} \label{eq:mainlemma}
 \ip{\{ f_S \}}{d} \leq  \tau(\epsilon,d) 
 \end{equation}
for a function $\tau(\epsilon,d) = \sqrt \epsilon \cdot 2^{O(d)}$ that we specify later.

We proceed by induction on $d$. 

\subsection*{The case $d=1$}
For $d=1$, we have two functions $f_\emptyset,f_{\{1\}}$ such that $\E_x f_{\{1\}}(x)=0$ and
we want an upper bound to $\ip{ \{ f_S \} } {1}$. We see that

\[ \ip { f_\emptyset,f_{\{1\}} } {1} = \( \E_x f_\emptyset (x) \) \cdot \( \E_y f_{\{1\}}(y) \)=0 \]

We have proved the base case of Equation \ref{eq:mainlemma} with $\tau(\epsilon,1) = 0$.

\subsubsection*{The inductive step}

Suppose now that, for $d\geq 1$, the lemma is true up to dimension $d$, and we want
to prove it for dimension $d+1$.

We have  $2^{d+1}$ functions $\{ f_S \}_{S\subseteq [d+1]}$ and we want to upper bound the Gowers
inner product

\begin{equation} \label{eq:ipdplusone} \ip { \{ f_S \} }{d+1} =  \E_{x,x_1,\ldots,x_{d+1}} \prod_{S\subseteq [d+1]} f_S\left(x+\sum_{i\in S} x_i\right)
\end{equation}

For every $x_1,\ldots,x_{d-1}$, define the four functions\footnote{The case $d=1$ is
somewhat degenerate: $[d-1]$ is the empty set, and $x_1,\ldots,x_{d-1}$ is an empty
sequence. So we simply have $A:= f_\emptyset$, $B:= f_{ \{ 1 \} }$, $C:= f_{ \{2 \} }$
and $D:= f_{ \{ 1,2 \} }$, with no subscripts.}
\[ A_{x_1,\ldots,x_{d-1}}(x) := \prod_{S\subseteq [d-1]} f_S \left ( x+ \sum_{i\in S} x_i \right) \]
\[ B_{x_1,\ldots,x_{d-1}}(x) := \prod_{S\subseteq [d-1]} f_{S\cup \{ d\} } \left ( x+ \sum_{i\in S} x_i \right) \]
\[ C_{x_1,\ldots,x_{d-1}}(x) := \prod_{S\subseteq [d-1]} f_{S \cup \{ d+1 \} } \left ( x+ \sum_{i\in S} x_i \right) \]
\[ D_{x_1,\ldots,x_{d-1}}(x) := \prod_{S\subseteq [d-1]} f_{S \cup \{ d,d+1 \} } \left ( x+ \sum_{i\in S} x_i \right) \]

with this notation, we can rewrite the expression~(\ref{eq:ipdplusone})  as

\[ \ip { \{ f_S \} }{d+1} = \E_{x,x_1,\ldots,x_{d+1}} A_{x_1,\ldots,x_{d-1}}(x) B_{x_1,\ldots,x_{d-1}}(x+x_d)
C_{x_1,\ldots,x_{d-1}}(x+x_{d+1})  D_{x_1,\ldots,x_{d-1}}(x+x_d + x_{d+1}) \]

and using the Fourier expansion and  simplifying,

\begin{eqnarray*} \ip { \{ f_S \} }{d+1} & = & \E_{x_1,\ldots,x_{d-1}} \sum_\alpha \hat A_{x_1,\ldots,x_{d-1}} (\alpha) \hat B_{x_1,\ldots,x_{d-1}} (\alpha) \hat C_{x_1,\ldots,x_{d-1}} (\alpha) 
\hat D_{x_1,\ldots,x_{d-1}} (\alpha) \\
& = & \E_{x_1,\ldots,x_{d-1}} \sum_{\alpha\neq \emptyset} \hat A_{x_1,\ldots,x_{d-1}} (\alpha) \hat B_{x_1,\ldots,x_{d-1}} (\alpha) \hat C_{x_1,\ldots,x_{d-1}} (\alpha) 
\hat D_{x_1,\ldots,x_{d-1}} (\alpha) \\ & &  + \E_{x_1,\ldots,x_{d-1}} \hat A_{x_1,\ldots,x_{d-1}} (\emptyset) \hat B_{x_1,\ldots,x_{d-1}} (\emptyset) \hat C_{x_1,\ldots,x_{d-1}} (\emptyset) 
\hat D_{x_1,\ldots,x_{d-1}} (\emptyset) \end{eqnarray*}

We bound the two terms separately.

For the first term, we have that, for every fixed $\bx = (x_1,\ldots,x_{d-1})$,

\begin{eqnarray*}  
& & \sum_{\alpha\neq \emptyset} \hat A_{\bx} (\alpha) \hat B_{\bx} (\alpha) \hat C_{\bx} (\alpha) 
\hat D_{\bx} (\alpha) \\
& \leq & 4 \max_{\alpha \neq \emptyset} \min\{ |\hat A_{\bx} (\alpha)|,| \hat B_{\bx} (\alpha)|,| \hat C_{\bx} (\alpha)| ,
|\hat D_{\bx} (\alpha)| \} \\
& \leq & 4 \sqrt{\max_i \min  \{\I_i (A_\bx) , \I_i (B_\bx), \I_i ( C_{\bx} ) , \I_i (D_{\bx}) \} }
\end{eqnarray*}

And we observe that

\begin{equation} \label{eq:inductivestep}
 \min  \{ \I_i (A_\bx) , \I_i (B_\bx),\I_i ( C_{\bx} ) , \I_i (D_{\bx}) \} \leq 2^{2d-2}\epsilon 
 \end{equation}
by using  \lemmaref{lm:influence-shift-general}. To verify this claim,\footnote{Again, the $d=1$
case is degenerate but easy to check:  Equation \ref{eq:inductivestep} reduces to
\[  \min  \{ \I_i (f_{\emptyset}) , \I_i (f_{ \{ 1 \} }),\I_i ( f_{ \{ 2 \} } ) , \I_i (f_{ \{ 1,2 \} }) \} \leq \epsilon\]
which is precisely our assumption that the 4-cross influence of the functions is at most $\epsilon$.} 
let $\delta$ be the minimum
in the above expression. Then $\I_i(A_\bx)\geq \delta$; recall $A_\bx$ is defined as a product
of $2^{d-1}$ functions of the form $f_S(x+ \sum_S x_i)$, $S\subseteq [d-1]$, where we think of the $x_i$ as constants,
and so for at least  one $S\subseteq [d-1]$ we have $\I_i (f_S) \geq \delta/2^{2d-2}$. We argue
similarly for $B_\bx$, $C_\bx$ and $D_\bx$, and we find three more functions for which
coordinate $i$ has influence at most $\delta/2^{2d-2}$, and these four functions are distinct.
Therefore, $\epsilon \geq 4\XI_i(\{ f_S \} ) \geq \delta/2^{2d-2}$. 

So we have that, for every fixed $\bx= (x_1,\ldots,x_{d-1})$

\[  \sum_{\alpha\neq \emptyset} \hat A_{\bx} (\alpha) \hat B_{\bx} (\alpha) \hat C_{\bx} (\alpha) 
\hat D_{\bx} (\alpha) \leq 4\sqrt {2^{2d-2}\epsilon} \]
and so, in particular
\[ \E_{\bx} \sum_{\alpha\neq \emptyset} \hat A_{\bx} (\alpha) \hat B_{\bx} (\alpha) \hat C_{\bx} (\alpha) 
\hat D_{\bx} (\alpha) \leq 4\sqrt {2^{2d-2}\epsilon} \]

The second term can be written as

\[ \E_{\bx} \hat A_{\bx} (\emptyset) \hat B_{\bx} (\emptyset) \hat C_{\bx} (\emptyset) 
\hat D_{\bx} (\emptyset) \]
\begin{equation} \label{eq:averaging} = \E_{x,\bx,y,z,w} A_{\bx}(x) B_{\bx}(x+y) C_{\bx}(x+z) D_{\bx} (x+w) \end{equation}
\begin{equation} \label{eq:after-averaging}
 \leq \E_{x,\bx,w} A_{\bx}(x) B_{\bx}(x+a) C_{\bx}(x+b) D_{\bx} (x+w) 
\end{equation}
where $a,b$ are values for $y,z$  that maximize the expectation in Equation \ref{eq:averaging}.

If we expand the Expression in (\ref{eq:after-averaging}), we get

\begin{equation} \label{eq:before:substitution}
\begin{array}{ll} 
\E_{x,x_1,\ldots,x_{d-1},w} \prod_{S\subseteq [d-1] } & f_S\left ( x + \sum_{i\in S} x_i \right)
f_{S\cup \{ d\} } \left ( x + a + \sum_{i\in S} x_i \right)\\
 & f_{S\cup \{ d+1\} } \left ( x + b + \sum_{i\in S} x_i \right)
f_{S\cup \{ d,d+1\} } \left ( x + w + \sum_{i\in S}x_i  \right)
\end{array}
\end{equation}
which we are going to re-write as the $d$-dimensional Gowers inner product
of a new set of functions, so that we can invoke the inductive hypothesis.
For every $S\subseteq [d-1]$, define
\[ g_S (x) := f_S(x) \cdot f_{S\cup \{ d \} } (x+a) \cdot f_{S \cup \{ d+1 \} } (x+b ) \]
and define
\[ g_{S\cup \{ d \} } (x) := f_{S\cup \{ d,d+1 \} } (x) \]

Expression (\ref{eq:before:substitution}) becomes
\[ \E_{x,x_1,\ldots,x_{d-1},w} \prod_{S\subseteq [d-1] } g_S\left( x + \sum_{\in S} x_i \right)\cdot
 g_{S \cup \{ d \} } \left(x+ \sum_{i \in S} x_i +w\right ) \]
\[ = \E_{x,x_1,\ldots,x_{d-1},x_d} \prod_{S\subseteq [d] } g_S\left( x + \sum_{\in S} x_i \right)
= \ip { \{ g_S \} } {d} \]
after the change of variable $w\to x_d$.

By definition, the function $g_{[d]}= f_{[d+1]}$ is balanced, and, by construction and
by \lemmaref{lm:influence-shift-general}, the 4-cross-influence of the functions $g_S$
is at most $9 \epsilon$. 
So we have $U^{d}(\{ g_S \}) \leq \tau\left(9 \epsilon,d\right)$.

We have thus solved the case of dimension $d+1$, with 
\[ \tau(\epsilon,d+1) = 2^{d+1}\cdot \sqrt \epsilon +  \tau\left( 9 \epsilon, d \right) \]
Together with the base case 
\[ \tau(\epsilon,1) = 0 \]
the recursion gives $\tau(\epsilon,d) = 2^{O(d)} \cdot \sqrt \epsilon$.
\end{proof}

\section{Generalizing to Complex-Valued Functions on Arbitrary Groups}
\label{sec:complex}

In this section we generalize our result to the setting of functions
$f: G_1 \times \cdots \times G_n \to \C$ where each $G_i$ is a finite abelian group.

We fix a group $G= G_1\times \cdots \times G_n$ for the rest of this section.

 We write
group operations in $G$ and in the groups $G_i$ additively. We will
recover the results of the previous sections when each $G_i$ is  
$\Z_2$. We denote an element
of $G$ as a tuple $\bg = (g_1,\ldots,g_n)$ where $g_i$ is an element of $G_i$.
The zero element of $G$ is the tuple $\bzero = (0,\ldots,0)$.

If $z= a+bi$ is a complex number, we define its  conjugate $\overline z:= a-bi$
and its absolute value
$|z|:= \sqrt{a^2 + b^2} = \sqrt{z \cdot \overline z}$.

\subsection{Complex Random Variables}

If $X$ is a random variable that takes on finitely many complex values,
and $\mu(x)$ is the probability that $X$ takes value $x$, then the
average of $X$ is the complex number

\[ \E [X] := \sum_{x} \mu(x) \cdot x \]

and the variance of $X$ is the {\em real} number

\[ \var[X] :=  \E[ | X- \E[X] |^2 ] \ . \]

\subsection{Fourier Analysis}

A function $\chi: G_0\to \C$ is a {\em character} of a group $G_0$
if $\chi(0) = 1$ and $\chi(a+b)= \chi(a)\chi(b)$ for every $a,b\in G_0$.
We also have $\chi(-a) = \overline{\chi(a)}$.

It is well known that a finite abelian group $G_0$ has precisely $|G_0|$
characters, and that there is an isomorphism between $G_0$ and the set
of characters that associates to each group element $g$ precisely
one character $\chi_g$ so that $\chi_0$ is the constant 1 function
and $\chi_g \chi_h = \chi_{g+h}$.

For each group $G_i$, let $\{ \chi^i_g \}_{g\in G_i}$ be the characters
of $G_i$, indexed according to the above isomorphism.

For each ${\bf g} = (g_1,\ldots,g_n)$ define the function $ \chi_{g_1,\ldots,g_n}: G \to \C$

\[ \chi_{g_1,\ldots,g_n} (x_1,\ldots,x_n) := \prod_{i\in [n]} \chi^i_{g_i} (x_i) \ . \]

Then the functions $\chi_{\bf g}$ are the characters of the group $G$.

Define the following inner product among functions $f,h:G \to \C$:

\[ \langle f, h\rangle := \E_x f(x) \overline{h(x)} \]

Then one can verify that the functions $\chi_{\bf g}$ are orthonormal
with respect to the inner product, and form a basis for the space
of functions $f:G \to \C$. Every such function can be written as

\[ f(x) = \sum_{{\bf g} \in G} \hat f(\bg) \chi_{\bg} (x) \]
where $\hat f(\bg) = \langle f, \chi_{\bg}\rangle$.

We again have Placherel's identity

\[ \E_x | f(x) |^2 = \sum_{\bg} |\hat f(\bg)|^2 \]

and we observe that

\[ \hat f(\bzero) = \E_x f(x) \]

\subsection{Influence of Variables}

The following definitions could be given in a much more general setting, but
the following will suffice for the purpose of this paper.

Let $\Sigma_1,\ldots,\Sigma_n$ be finite and $f:\Sigma_1\times \cdots \times \Sigma_n \to \C$ 
be a function, then the {\em influence of $i$} for $f$
is defined as

\[ \I_i(f) := \E_{x_1,\ldots,x_{i-1},x_{x+1},\ldots,x_n} \left[ \var_{x_i} [f(x_1,\ldots,x_n)] \right] \]

where the variables $x_i$ are mutually independent, and each $x_i$ is uniformly distributed over $\Sigma_i$. 
The reader should verify that
if each $\Sigma_i = \B$ and if $f$ takes only real values, then we recover the definition we gave earlier.
Note that the influence is always a non-negative real number.

Returning to our setting of functions $f:G_1\times \cdots \times G_n\to \C$, we have the following fact:

\[  \I_i(f) = \sum_{(g_1,\ldots,g_n) : g_i \neq 0} | \hat f(g_1,\ldots,g_n) |^2 \]

and, in particular,

\[ \max_i \I_i (f) \geq \max_{\bg \neq \bzero} | \hat f(\bg) |^2 \]

\subsection{Influence of Products of Functions}

The following result is the only real difficulty in generalizing our results
from previous sections. (The rest just follows the same proofs with some
changes in notation.) 

\begin{lemma} \label{lm:influence-complex} \label{lm:influence-shift-complex}
Let $f,g:\Sigma_1\times \cdots \times \Sigma_n \to \C$ be functions such that $|f(x)|\leq 1$ and $|g(x)|\leq 1$
for every $x$.

Then, for every $i\in [n]$,
\[ \I_i (fg) \leq 3 \cdot (\I_i(f) + \I_i(g) ) \]
\end{lemma}

To prove the lemma, it is enough to prove the following bound

\begin{lemma}
Let $f,g:\Sigma \to \C$ be functions such that $|f(x)|\leq 1$ and $|g(x)|\leq 1$
for every $x$.

Then,
\[ \var_x [f(x)g(x)] \leq 3 \cdot (\var_x[f(x)]  + \var_x[g(x)] ) \]
\end{lemma}

\begin{proof}
Recall that $\var_x [ f(x)g(x) ]  :=  \E_x \left| f(x)g(x) - \E_y [f(y)g(y)] \right|^2$.
For every $x$,

\begin{eqnarray*}
\left| f(x)g(x) - \E_y [f(y)g(y)] \right| 
& = &  \left| f(x)g(x) - \E_y[f(y)]g(x) 
 +\E_y[f(y)]g(x)- \E_y [f(y)g(y)] \right| \\
& = &  \left| g(x) \cdot \left( f(x) - \E_y[f(y)] \right) + \E_y\left[f(y)\cdot (g(x)-g(y)) \right]\right|\\
& \leq & | g(x) | \cdot \left| f(x) - \E_y[f(y)] \right| + \E_y\left[|f(y)| | g(x)-g(y) | \right]\\
& \leq & \left| f(x) - \E_y[f(y)] \right| + \E_y\left[|f(y)| \left| g(x)-g(y) \right| \right]
\end{eqnarray*}

We can now compute the variance of $fg$:

\begin{eqnarray*}
\var_x [ f(x)g(x) ] & := & \E_x \left| f(x)g(x) - \E_y[f(y)g(y)] \right|^2\\
& \leq & \E_x \left| \left| f(x) - \E_y[f(y)] \right| + \E_y\left[|f(y)| \left| g(x)-g(y) \right| \right] \right|^2\\
& \leq & 2\E_x\left[ \left| f(x) - \E_y[f(y)] \right|^2 + \left| \E_y \left[|f(y)| \left| g(x)-g(y) \right| \right] \right|^2 \right]\\
& = & 2\E_x \left| f(x) - \E_y[f(y)] \right|^2  + 2\E_x\E_y \left[ |f(y)|^2 |g(x)-g(y)|^2 \right]\\
& \leq & 2\var_x [f(x)] + 4\var_x[g(x)]
\end{eqnarray*}

Where the last step follows from the fact that

\[ \E_{x,y} |g(x)-g(y)|^2 = \E_{x,y} \left| g(x) - \E_z[g(z)] + \E_z[g(z)] - g(y) \right|^2
\leq 2 \E_x \left| g(x) - \E_z[g(z)] \right|^2 + 2\E_y \left| \E_z[g(z)] - g(y) \right|^2 \]

Similarly, we could prove
\[ \var_x [ f(x)g(x) ] \leq 4\var_x [f(x)] + 2\var_x[g(x)] \]
And the average of the two bonds gives us the desired result.
\end{proof}

A simple induction shows that
 
\[ \I_i (f_1 \cdots f_k) \leq 3\cdot k^{\log_2 3} \cdot \sum_j \I_i(f_j) \]

\subsection{Gowers Uniformity}

The definition of $U^d$ for complex-valued functions is as follows:
for $d\geq 1$ and $f:G\to \C$,

\[ U^d(f) := \E_{x,x_1,\ldots,x_d} \left(
 \prod_{S\subseteq [d], |S|~{\rm even}} f\left(x+ \sum_{i\in S} x_i\right) \right) \cdot
 \left(
 \prod_{S\subseteq [d], |S|~{\rm odd}} \overline{f\left(x+ \sum_{i\in S} x_i\right)} \right)
\]

The inductive definition of $U^d$ is perhaps simpler.
For a function $f:G \to \C$ and elements $x_1,\ldots,x_d \in G$,
define $f_{x_1,\ldots,x_d} : G \to \C$ inductively as follows:

\[ f_{x_1,\ldots,x_d} (x) := f_{x_1,\ldots,x_{d-1}} (x) \overline{f_{x_1,\ldots,x_{d-1}} (x+x_d)} \]

and then define

\[ U^d (f) :=  \E_{x,x_1,\ldots,x_d} f_{x_1,\ldots,x_d} (x) \]
It is possible to show that $U^d(f)$ is always a non-negative real number, because we have

\begin{eqnarray*}
 U^d(f) & = & \E_{x,x_1,\ldots,x_d} f_{x_1,\ldots,x_{d-1}} (x)\overline {f_{x_1,\ldots,x_{d-1}} (x+x_d)} \\
& = & \E_{x_1,\ldots,x_{d-1}} \left( \E_x f_{x_1,\ldots,x_{d-1}} (x) \right) 
\left(\overline{ \E_x f_{x_1,\ldots,x_{d-1}} (x) } \right) \\
& = & \E_{x_1,\ldots,x_{d-1}} \left| \E_x f_{x_1,\ldots,x_{d-1}} (x) \right|^2
\end{eqnarray*}

Explicit formulas for the case $d=1,2$ are:

\[ U^1(f) = \E_{x,y} f(x) \overline { f(x+y) } = (\E_{x} f(x)) \cdot \overline {\E_x f(x))} = | \E_x f(x) |^2 \]

and

\[ U^2(f) = \E_{x,y,z} f(x) \overline {f(x+y)} \overline {f(x+z)} f(x+y+z)
= \sum_{\bg} \hat f(\bg) \overline{  \hat f(\bg)} \overline{  \hat f(\bg)} \hat f(\bg)  = \sum_{\bg} |\hat f(\bg) |^4 \]

\subsection{Gowers Uniformity and Influence}

\begin{theorem}
Let $f:G_1\times \cdots \times G_n \to \C$ be a function such that $|f(x)| \leq 1$ for every $x$. Then, for every $d\geq 1$,

\[ U^d(f) \leq U^1(f) +2^{O(d)}  \max_i \I_i (f) \]
\end{theorem}

\begin{proof} The case $d=1$ is trivial. For $d\geq 2$, it suffices to 
prove $$U^d(f) \leq U^{d-1}(f) + 2^{cd} \max_i \I_i (f)$$ for an absolute constant $c$. Let
$I:= \max_i \I_i(f)$.

We follow the proof of \lemmaref{lm:boolean-influence-one} and write

\[ U^d(f) =  \E_{x_1,\ldots,x_{d-2}}\left [ \E_{x,y,z} f_{x_1,\ldots,x_{d-2}}(x)
\overline{f_{x_1,\ldots,x_{d-2}}(x+y)} \overline {f_{x_1,\ldots,x_{d-2}}(x+z) }f_{x_1,\ldots,x_{d-2}}(x+y+z) \right] \]
 \[ = \E_{x_1,\ldots,x_{d-2}}\sum_{\bg} | \hat f_{x_1,\ldots,x_{d-2}}(\bg ) |^4\]
\[ = \E_{x_1,\ldots,x_{d-2}}\sum_{\bg \neq \bzero} | \hat f_{x_1,\ldots,x_{d-2}}(\bg ) |^4
+ \E_{x_1,\ldots,x_{d-2}} | \hat f_{x_1,\ldots,x_{d-2}}( \bzero) |^4 \]

We separately bound the two terms in the last expression. 

For every $a_1,\ldots,a_{d-2} \in G$, and every $\bg= (g_1,\ldots,g_n) \neq \bzero$, using 
\lemmaref{lm:influence-complex}, we get that, for every $i: g_i \neq 0$,
$\I_i(f_{a_1,\ldots,a_{d-2}}) \leq 2^{cd} \I_i (f) \leq 2^{cd} I$ for an absolute constant $c$,
and so we have that, for every $a_1,\ldots,a_{d-2} \in G$,

\[ \sum_{\bg\neq \bzero} | \hat f_{a_1,\ldots,a_{d-2}}(\bg) |^4
\leq \max_\bg | \hat f_{a_1,\ldots,a_{d-2}}(\bg ) |^2 \leq  \max_i \I_i (f_{a_1,\ldots,a_{d-2}})
\leq 2^{cd} I \]
and so
\[ \E_{x_1,\ldots,x_{d-2}}\sum_{\bg\neq \bzero} | \hat f_{x_1,\ldots,x_{d-2}}(\bg ) |^4 \leq 2^{cd} I \]

Regarding the other term,
\[ \E_{x_1,\ldots,x_{d-2}} | \hat f_{x_1,\ldots,x_{d-2}}(\bzero) |^4 \leq 
\E_{x_1,\ldots,x_{d-2}} | \hat f_{x_1,\ldots,x_{d-2}}(\bzero) |^2
 = \E_{x_1,\ldots,x_{d-2}} \left| \E_x f_{x_1,\ldots,x_{d-2}}(x )\right |^2 = U^{d-1}(f) 
\]
\end{proof}

\subsection{Gowers Inner Product}

Let $\{ f_S \}_{S\subseteq [d]}$ be a collection of functions $f_S: G \to \C$. Then
their Gowers inner product is the complex number

\begin{equation}  \ip { \{ f_S \} }{d} :=  \E_{x,x_1,\ldots,x_{d}} \left[ 
\prod_{S\subseteq [d]: |S| ~ \rm even} f_S\left(x+\sum_{i\in S} x_i\right) \cdot
\prod_{S\subseteq [d]: |S| ~ \rm odd} \overline{f_S\left(x+\sum_{i\in S} x_i\right)} \right]
\end{equation}

\subsection{Gowers Inner Product and Cross-Influence}

We generalize \lemmaref{lm:main} to the case of products of arbitrary groups.

\begin{lemma} 
For $d\geq 1$, let $\{ f_S \}_{S\subseteq [d]}$ be a collection
of functions $f_S: G_1\times \cdots \times G_n \to \C$ such that
\begin{itemize}
\item $|f_S(x)|\leq 1$ for every $S$ and every $x$;
\item $f_{[d]}$ is balanced, that is, $\E_x f_{[d]}(x) = 0$;
\item $4\XI_i (\{ f_S \}) \leq \epsilon$ for every $i$.
\end{itemize}

Then 
\[ | \ip{\{ f_S \} }{d} | \leq \sqrt{\epsilon} \cdot 2^{O(d)} \]
\end{lemma}

We give an outline of the changes  needed to adapt the proof of \lemmaref{lm:main}.

First, we need the following bound, whose proof is identical to the proof
of \lemmaref{lm:fourfunctions}.

\begin{lemma}\label{lm:fourfunctions-c}
Let $f_1,f_2,f_3,f_4 : G \to \C$ be functions such that $|f_i(x)|\leq 1$ for every $i$ and every $x$.
Then

\[ \sum_{\bg} |\hat f_1(\bg) |\cdot |\hat f_2(\bg) | \cdot |\hat f_3(\bg) | \cdot |\hat f_4(\bg) |
\leq 4\max_\bg \min\{ |\hat f_1(\bg) |, |\hat f_2(\bg) | , |\hat f_3(\bg) | , |\hat f_4(\bg) | \} 
\]
\end{lemma}

As in the proof \lemmaref{lm:main}, we proceed by induction on $d$ and
prove that, under the hypothesis of the lemma, $\ip{\{ f_S \} }{d} | \leq \tau(\epsilon,d)$
for a  function $\tau$ that satisfies $\tau(\epsilon,d) = \sqrt \epsilon 2^{O(d)}$. For $d=1$,
We have $\ip{ f_\emptyset, f_ { \{ 1 \} } } {1} = \( \E_x f_\emptyset (x) \) \cdot \(
\E_y \overline{f_{ \{ 1 \} } (y)} \) = 0$, and so we have the base case with $\tau(\epsilon,1)=0$.

For the inductive step, we consider the $(d+1)$-dimensional Gower inner product

\begin{equation}  \ip { \{ f_S \} }{d+1} =  \E_{x,x_1,\ldots,x_{d+1}} \left[ 
\prod_{S\subseteq [d+1]: |S| ~ \rm even} f_S\left(x+\sum_{i\in S} x_i\right) \cdot
\prod_{S\subseteq [d+1]: |S| ~ \rm odd} \overline{f_S\left(x+\sum_{i\in S} x_i\right)} \right]
\end{equation}

It is easier to rewrite it as

\begin{equation}  \label{eq:ipdplusone-c}\ip { \{ f_S \} }{d+1} =  \E_{x,x_1,\ldots,x_{d+1}} \left[ 
\prod_{S\subseteq [d+1]: |S|} F_S\left(x+\sum_{i\in S} x_i\right) \right]
\end{equation}
where we define $F_S:= f_S$ if $|S|$ is even and $F_S:= \overline {f_S}$ if $|S|$ is odd.

For every $x_1,\ldots,x_{d-1}$, define the four functions
\[ A_{x_1,\ldots,x_{d-1}}(x) := \prod_{S\subseteq [d-1]} F_S \left ( x+ \sum_{i\in S} x_i \right) \]
\[ B_{x_1,\ldots,x_{d-1}}(x) := \overline{\prod_{S\subseteq [d-1]} F_{S\cup \{ d\} } \left ( x+ \sum_{i\in S} x_i \right)} \]
\[ C_{x_1,\ldots,x_{d-1}}(x) := \overline{\prod_{S\subseteq [d-1]} F_{S \cup \{ d+1 \} } \left ( x+ \sum_{i\in S} x_i \right)} \]
\[ D_{x_1,\ldots,x_{d-1}}(x) := \prod_{S\subseteq [d-1]} F_{S \cup \{ d,d+1 \} } \left ( x+ \sum_{i\in S} x_i \right) \]

with this notation, we can rewrite the expression~(\ref{eq:ipdplusone-c})  as

\[ \ip { \{ f_S \} }{d+1} = \E_{x,x_1,\ldots,x_{d+1}} A_{x_1,\ldots,x_{d-1}}(x) \overline{B_{x_1,\ldots,x_{d-1}}(x+x_d)}
\overline{C_{x_1,\ldots,x_{d-1}}(x+x_{d+1})}  D_{x_1,\ldots,x_{d-1}}(x+x_d + x_{d+1}) \]

and using the Fourier expansion and  simplifying,

\begin{eqnarray*} | \ip { \{ f_S \} }{d+1} | & = & 
\left|\E_{\bx=(x_1,\ldots,x_{d-1})} \sum_\bg 
 \hat A_{\bx} (\bg)
 \overline{ \hat B_{\bx} (\bg) }
 \overline{ \hat C_{\bx} (\bg) }
 \hat D_{\bx} (\bg) \right|\\
& \leq & \left | \E_{\bx= (x_1,\ldots,x_{d-1})} \sum_{\bg\neq (0,\ldots,0)} 
 \hat A_{\bx} (\alpha)
 \overline{ \hat B_{\bx} (\bg) }
 \overline{ \hat C_{\bx} (\bg) }
 \hat D_{\bx} (\bg) \right| \\ 
& &  + \left|  \E_{\bx = (x_1,\ldots,x_{d-1})}  \hat A_{\bx} (0,\ldots,0)
 \overline{ \hat B_{\bx} (0,\ldots,0) }
 \overline{ \hat C_{\bx} (0,\ldots,0) }
 \hat D_{\bx} (0,\ldots,0) \right| \end{eqnarray*}

We bound the two terms separately.

For the first term, we have that, for every fixed $\bx = (x_1,\ldots,x_{d-1})$,

\begin{eqnarray*}  
& & \left| \sum_{\bg\neq (0,\ldots,0)} \hat A_{\bx} (\bg) \overline {\hat B_{\bx} (\bg) \hat C_{\bx} (\bg)} 
\hat D_{\bx} (\bg)\right| \\
& \leq & \sum_{\bg\neq (0,\ldots,0)} | \hat A_{\bx} (\bg) | \cdot | \overline {\hat B_{\bx} (\bg)}| \cdot
|\overline{ \hat C_{\bx} (\bg)} | \cdot | 
\hat D_{\bx} (\bg) |\\
& = & \sum_{\bg\neq (0,\ldots,0)} | \hat A_{\bx} (\bg) | \cdot |  {\hat B_{\bx} (\bg)}| \cdot
|{ \hat C_{\bx} (\bg)} | \cdot | 
\hat D_{\bx} (\bg) |\\
& \leq &  4 \sqrt { \max_i \min  \{ \I_i ( A_{\bx} ),\I_i ( B_{\bx} ),\I_i ( C_{\bx} ) , \I_i (D_{\bx}) \} }\\
\end{eqnarray*}

And we observe that, for every $\bx=(x_1,\ldots,x_{d-1})$ and every  $i$,

\[ \min  \{ \I_i ( A_{\bx} ),\I_i ( B_{\bx} ),\I_i ( C_{\bx} ) , \I_i (D_{\bx}) \} \leq 2^{cd}\epsilon \]
for some absolute constant $c$.

The second term can be written as

\[ \left| \E_{\bx} \hat A_{\bx} (0,\ldots,0) \overline{\hat B_{\bx} (0,\ldots,0) \hat C_{\bx} (0,\ldots,0) } 
\hat D_{\bx} (0,\ldots,0) \right|  \]
\begin{equation} \label{eq:averaging-c} = \left| \E_{x,\bx,y,z,w} A_{\bx}(x)\overline {B_{\bx}(x+y) C_{\bx}(x+z) }D_{\bx} (x+w) \right| \end{equation}
\[ \leq \E_{y,z} \left | \E_{x,\bx,w} A_{\bx}(x)\overline {B_{\bx}(x+y) C_{\bx}(x+z) }D_{\bx} (x+w) \right| \]
\begin{equation} \label{eq:after-averaging-c}
 \leq \left| \E_{x,\bx,w} A_{\bx}(x) \overline{B_{\bx}(x+a) C_{\bx}(x+b)} D_{\bx} (x+w) \right|
\end{equation}
where $a,b$ are values for $y,z$  that maximize the expectation in Equation \ref{eq:averaging-c}.

If we expand the Expression in (\ref{eq:after-averaging-c}), we get

\begin{equation} \label{eq:before:substitution-c}
\begin{array}{ll} 
\E_{x,x_1,\ldots,x_{d-1},w} \prod_{S\subseteq [d-1] } & F_S\left ( x + \sum_{i\in S} x_i \right)
F_{S\cup \{ d\} } \left ( x + a + \sum_{i\in S} x_i \right)\\
 & F_{S\cup \{ d+1\} } \left ( x + b + \sum_{i\in S} x_i \right)
F_{S\cup \{ d,d+1\} } \left ( x + w + \sum_{i\in S}x_i  \right)
\end{array}
\end{equation}
which we are going to re-write as the $d$-dimensional Gowers inner product
of a new set of functions, so that we can invoke the inductive hypothesis.
For every $S\subseteq [d-1]$, if $S$ is even, define
\[ g_S (x) := F_S(x) \cdot F_{S\cup \{ d \} } (x+a) \cdot F_{S \cup \{ d+1 \} } (x+b ) \]
and define
\[ g_{S\cup \{ d \} } (x) := \overline{F_{S\cup \{ d,d+1 \} } (x)} \]
If $S$ is odd, define
\[ g_S (x) := \overline{ F_S(x) \cdot F_{S\cup \{ d \} } (x+a) \cdot F_{S \cup \{ d+1 \} } (x+b ) } \]
and 
\[ g_{S\cup \{ d \} } (x) := {F_{S\cup \{ d,d+1 \} } (x)} \]

Expression (\ref{eq:before:substitution-c}) becomes
$\ip { \{ g_S \} } {d}$
after the change of variable $w\to x_d$.

By definition, the function $g_{[d]}$ is either $f_{[d+1]}$ of $\overline{f_{[d+1]} }$ and,
in either case, $g_{[d]}$ is balanced. By construction and
by \lemmaref{lm:influence-shift-complex}, the cross-influence of the functions $g_S$
is at most $c' \epsilon$ for an absolute constant $c'$. 

This proves the inductive step with $\tau(\epsilon, d+1) = 4\cdot \sqrt{2^{cd} \epsilon} + \tau(c'\epsilon,d-1)$
where $c,c'$ are absolute constants.
Together with $\tau(\epsilon,1) = 0$ we have $\tau(\epsilon,d) = 2^{O(d)} \sqrt \epsilon$
as desired.

\section{A Tight Analysis of Linearity Testing}
\label{sec:lowerbound}

Consider the following promise problem. Given a function $f:\B^n \to \H$ and a small $\epsilon>0$, 
we want to distinguish the two cases

\begin{enumerate}

\item $f$ is linear;
\item $U^d(f) \leq \epsilon$.
\end{enumerate}
We refer to such a test as a ``relaxed linearity test of degree $(d-1)$.'' 
As usual, we say that a test has completeness $c$ and soundness $s$ if the test acceptes with
probability $\geq c$ in case (1) and with probability $\leq s+\epsilon'$ in case
(2), where $\epsilon' \to 0$ when $\epsilon \to 0$. If a test makes $q$
queries and has soundness $s$, then its {\em amortized query complexity}
is $\bar q = q/\log s^{-1}$. 

For $d=2$, this problem is the linearity testing problem. For $d=3$, the
only functions such that $U^3(f)\geq \epsilon$ are functions that are correlated
with degree-2 polynomials \cite{S05}, and so the test is required to
distinguish linear functions from functions that are far from being quadratic.
For $d\geq 4$, it is conjectured that the only functions with $U^d(f)\geq \epsilon$
are those that are correlated with a degree-$(d-1)$ polynomial, and, if so, such
a test distinguishes linear functions from functions that are far from low-degree
polynomials. By our results, such a test also distinguishes linear functions
from functions where all variables have low influence.

We give a tight analysis of the error probability of such tests for a given
number of queries.

\subsection{The Linear Gowers Inner Product}

For the sake of our analysis of the Hypergraph Test, it is convenient
to study expressions of the following form. Let $\{ f_S \}_{S\subseteq [d]}$
be a collection of $2^d$ functions $f_S:\B^n \to [-1,1]$, $d\geq 1$, and define their
{\em linear} Gowers Inner Product as

\[ \lip{ \{ f_S \} }{d} := \E_{x_1,\ldots,x_d} \prod_{S\subseteq [d]} f \( \sum_{i\in S} x_i \) \]
As usual, an empty sum is zero. For example:

\begin{eqnarray*}
\lip { f_\emptyset , f_{ \{ 1 \} } } {1} & := &  \E_{x} f_\emptyset (\bzero) f_{ \{ 1 \} } (x) \\
\lip { f_\emptyset , f_{ \{ 1 \} }, f_{ \{ 2 \} } , f_{ \{ 1,2 \} } } {2} 
& := & \E_{x,y} f_\emptyset (\bzero) f_{ \{ 1 \} } (x) f_{ \{ 2 \} } (y)f_{ \{ 1,2 \} } (x+y)\\
\end{eqnarray*} 

Where $\bzero= ( 0,\ldots,0)$ is the all-zero vector of $\B^n$.
We call it a {\em linear} inner product because, for functions $f_S:\B^n \to \R$, the Gowers inner
product is defined by picking at random an {\em affine} subspace of $\B^n$ of dimension $d$, and
then taking the product of the functions on all points of the subspace. In the above expression, we
do something similar but on a {\em linear} subspace.

We prove that if the linear Gowers inner product of a collection of functions is large,
then the regular Gowers inner product of a related collection of functions must also be large.

\begin{lemma}\label{lm:lingowers}
Let $f_S:\B^n \to [-1,1]$ be functions, $S\subseteq [d]$, and define
the collection $\{ g_T \}_{T\subseteq [d]}$ as $g_T = f_{T\cup \{ d\}}$. Then

\[ |\lip { \{ f_S \} }{d}| \leq \sqrt{ \ip{ \{ g_S \} } {d}} \]
\end{lemma}

\begin{proof}
\begin{eqnarray*}
|\lip { \{ f_S \} }{d}| & := &  \left|\E_{x_1,\ldots,x_d} \prod_{S\subseteq [d]} f_S \( \sum_{i\in S} x_i \) \right|\\
& \leq & \sqrt{\E_{x_1,\ldots,x_{d-1}} \( \prod_{S\subseteq [d],  d\not\in S} f_S \( \sum_{i\in S} x_i \)\)^2}
\times\\
& & \sqrt{\E_{x_1,\ldots,x_{d-1}} \( \E_{x_d} \prod_{S\subseteq [d],  d\in S} f_S \( \sum_{i\in S} x_i \)\)^2}\\
& \leq & \sqrt{\E_{x_1,\ldots,x_{d-1}} \( \E_{x_d} \prod_{S\subseteq [d],  d\in S} f_S \( \sum_{i\in S} x_i \)\)^2}\\
& = & \sqrt{\E_{x_1,\ldots,x_{d-1}}  \( 
\E_{x} \prod_{S\subseteq [d],  d\in S} f_S \( x+\sum_{i\in S-\{ d \}} x_i \) \) \cdot
\( \E_{y} \prod_{S\subseteq [d],  d\in S } f_S \( y+\sum_{i\in S-\{ d \}} x_i \) \)}\\
& = & \sqrt{ \E_{x,y,x_1,\ldots,x_{d-1}} \prod_{T\subseteq [d-1]} f_{T\cup \{ d\} } \( x+ \sum_{i\in T} x_i \)
f_{T \cup \{ d \} } \( y+ \sum_{i\in T} x_i \)}\\
& \leq & \sqrt{\ip{ \{ g_T \}}{d}} 
\end{eqnarray*}
after the change of variable $y\to x+x_d$.
\end{proof}

\begin{lemma}\label{lm:lip}
Let $f_S :\B^n \to\J$ be functions, $S\subseteq [d]$, such that
$\lip{ \{ f_S \} } d \geq \epsilon$ and $\E_x f_{[d]} (x) = 0$. 

Then there is a variable $i$ such that $\XI_i (\{ f_S \} ) \geq \epsilon^4/2^{O(n)}$.
\end{lemma}

\begin{proof}
Define $g_T:= f_{T\cup \{ d \} }$. Then we have
\begin{enumerate}
\item $\E_x g_{[d]} (x) = \E_x f_{[d]} (x) = 0$,
\item $\ip{ \{ g_T \} }d \geq ( \lip{ \{ f_S \} } d )^2 \geq \epsilon^2$.
\end{enumerate}
{From} \lemmaref{lm:main} we derive that there is a variable $i$
such that $4\XI_i ( \{ g_T \} ) \geq \epsilon^4/2^{O(d)}$. Each function $f_S$
occurs at most twice in the collection $\{g_T\}$, and so it must be
$\XI_i ( \{ f_S \} ) \geq 4\XI_i ( \{ g_T \} ) \geq \epsilon^4/2^{O(d)}$.
\end{proof}

\subsection{Positive Results on Relaxed Linearity Testing}
\label{sec:positive-linear}

Given a hypergraph $H = ([k],E)$, we can define a relaxed linearity test
associated with $H$ by 

\noindent\fbox{\begin{minipage}{2in}
\begin{program}
$H$-Test\+\\
 choose $x,x^{1},\ldots,x^{k}$ uniformly at random in $\B^n$\\
 \ACCEPT\ if and only if\+\\ $\forall e \in E: \ \prod_{i\in e}
 f\(x^i\) = f\(\sum_{i \in e} x^{i}\)$ 
\end{program}\end{minipage}}

Then we have the following result.

\begin{theorem} \label{th:relaxedlinear}
Let $d\geq 2$ and let $H=([k],E)$ be an hypergraph
such that each edge  of $H$ contains at most $d$ vertices.

Them the $H$-Test is a degree-$(d-1)$ relaxed linearity test of completeness 1 and
soundness at most $1/2^{|E|}$.
\end{theorem}

We remark that this result was first proved in \cite{S05}, using a 
different approach.

For the proof of \autoref{th:relaxedlinear} and of results in the next section,
it will be convenient to use the following ``Vazirani XOR Lemma,'' whose proof is immediate.

\begin{lemma} \label{lm:vazirani}
Let $X_1,\ldots,X_m$ be random variables taking values in $\H$. Then

\[ \pr [ X_1=1 \wedge X_2 = 1 \wedge \cdots \wedge X_m = 1 ] = \frac {1}{2^m}\sum_{S\subseteq [m]} \prod_{i\in S} X_i \]
\end{lemma}

Using \lemmaref{lm:vazirani}, we see that the probability that the $H$-Test accepts
a function $f$ is equal to

\[ \frac 1 {2^{|E|} } \sum_{E' \subseteq E} \prod_{e\in E} \( \prod_{i\in e} f(x^i) \) \cdot f\( \sum_{i\in e} x^i \) \]
\[ =  \frac 1 {2^{|E|} } +  \frac 1 {2^{|E|} }\sum_{E' \subseteq E, E'\neq \emptyset} \prod_{e\in E} \( \prod_{i\in e} f(x^i) \) \cdot f\( \sum_{i\in e} x^i \)\]

We will also need two results from \cite{Gowers01,GT:primes}. The first one is  \cite[Lemma 3.8]{Gowers01},
and it states that for every collection $\{ f_S \}$ of functions:
$$
\Big | \ip{ \{ f_S \} } {d}  \Big | \leq  
 \prod_S \( U^d\( f_S\) \)^{1/2^d}
$$
The other is \cite[5.7]{GT:primes} and states that for every $f$ and $d\geq 2$,
\[ U^{d-1}(f)  \leq \sqrt{ U^d (f) }\]

We now proceed with the proof of \autoref{th:relaxedlinear}.

\begin{proof}[Of \autoref{th:relaxedlinear}]
It's clear that a linear function is accepted with probability 1.

If 
the $H$-Test accepts with $f$ probability more than
$1/2^{|E|} + \epsilon$, then there is a non-empty $E'\subseteq E$ such that
\begin{equation} \label{eq:linear-start}  \prod_{e\in E'} \( \prod_{i\in e} f(x^i) \) \cdot f\( \sum_{i\in e} x^i \)
\geq \epsilon \end{equation}

Let $d'$ be the size of the largest edge in $E'$, and, without loss of
generality, assume that the edge $(1,\ldots,d')$ is in $E'$. Fix 
the variables $x^{d'+1},\ldots,x^k$ to values that maximize (\ref{eq:linear-start}).
Then, (\ref{eq:linear-start}) becomes

\[ \prod_{S\subseteq [d']}  f_S \( \sum_{i \in S} x_i \) \geq \epsilon \]
where $f_\emptyset$ is the constant function equal to the
product of the terms of (\ref{eq:linear-start}) that depend only on $x^{d'+1},\ldots,x^k$, $f_{ \{ 1 \} }(x^1 )$
is the product of the  terms of (\ref{eq:linear-start}) that depend only on $x^1$ and $x^{d'+1},\ldots,x^k$,
and so on. In particular, $f_{[d']} ( x^1 + \cdots + x^{d'} ) = f( x^1 + \cdots + x^{d'} ) $.

We thus have

\[ \lip { \{ f_S \} } {d'} \geq \epsilon \]
By  \lemmaref{lm:lingowers}, there are functions $\{ g_S \}$
such that $g_{[d']} = f$ and $\ip { \{ g_S \} } {d'} \geq \epsilon^2$.

Since all the functions involved are boolean, and thus their
uniformity norms are at most $1$, we conclude
$$
 \epsilon^2 \leq \Big| \ip{ \{ g_S \} } {d'}  \Big |  \leq \min_S
|U^{d'}\(g_S\)|^{1/2^{d'}} \le |U^{d'}\(f\)|^{1/2^{d'}}
$$
And, finally, $U^d (f) \geq \( U^{d'} (f) \)^{2^{d-d'}} \geq \epsilon^{2^{d+1}}$.
\end{proof}

In particular, taking $H$ to be the complete hypergraph on $k$
vertices with at most $d$ vertices per edge, 
leads to a hypergraph linearity test with $\sum_{i=1}^d {k \choose i}$
queries, and 
soundness $\frac{1}{2^{\sum_{i=2}^d {k \choose i}}}$. The amortized query
complexity of this test, for $k \gg d$, is $\bar{q} \le 1 + 
O\(\frac{1}{q^{(d-1)/d}}\)$.

Next, we show that no linearity test can do better. 

We do it in two steps. The first is to show that 
the amortized query complexity of any hypergraph test 
cannot be better than $1 + \Omega\(\frac{1}{q^{(d-1)/d}}\)$ for
$q$ queries. 

The second step builds on the first, and 
shows that {\it any} non-adaptive linearity test with perfect
completeness cannot do better under this promise, namely it will have
amortized query 
complexity of at least $1 + \Omega\(\frac{1}{q^{(d-1)/d}}\)$ for
$q$ queries. 

\subsection{Lower Bound for the $H$-Test}

We  prove a lower bound for the $H$-Test by describing an explicit function 
$f:~\B^n\rightarrow
\H$, which has small $d$-th uniformity norm, and for any
hypergraph $H$ the 
acceptance probability of the $H$-Test  
on $f$ is large. 
\begin{proposition}
\label{pro:hypertest} Fix $d\geq 2$ and fix an hypergraph $H = ([k],E)$
such that all edges have at most $d$ vertices.
For $n$ a multiple of $d$, let   $f: \B^n \to \H$ be defined as 
$$ 
f(x_1,\ldots,x_n) := (-1)^{x_1 x_2 \cdots x_d + x_{d+1} \cdots
  x_{2d} + \cdots}
$$ 
Then $U^d(f) \le 2^{-\Omega(n)}$, and 
the $H$-Test  accepts $f$ with probability at
least 
$$
\max\left\{1/2^{|E|}, 2^{-\sum_{i=2}^d {k \choose i}}\right\}
$$
\end{proposition}
The lower bound on amortized query complexity of hypergraph tests
follows immediately, if we recall that the number of queries is $|E| +
k$, and consider the two cases: $|E| \le \sum_{i=2}^d {k \choose i}$,
  or $|E| > \sum_{i=2}^d {k \choose i}$.

\begin{proof}
We first show $U^d(f)$ to be small. 

Some additional notation: 
Let $g_i(x) = x_{i\cdot d + 1}
  x_{i\cdot d + 2} \cdots x_{(i+1) \cdot d}$, and let $g(x) = x_1 x_2
  \cdots x_d +  x_{d+1} \cdots x_{2d} + \cdots = \sum_{i=0}^{n/d-1}
  g_i$. Let $f_i = (-1)^{g_i}$. 

Then $f = \prod_{i=0}^{n/d-1} f_i$ and $U^d(f) = \prod_{i=0}^{n/d-1}
U^d(f_i) = \(U^d(f_0)\)^{n/d}$. Thus it remains to show that
$U^d(f_0)$ is bounded away from $1$. 
$$
U^d(f_0) = U^d\((-1)^{g_0}\) = \E_{x,x^1,...,x^d} (-1)^{\sum_{S
    \subseteq [d]} g_0\(x + \sum_{i\in S} x^i\)} 
$$
We may assume that the variables $x^i$ live in $\B^d$, and then $g_0$
is just the $AND$ function. Therefore $\sum_{S
    \subseteq [d]} g_0\(x + \sum_{i\in S} x^i\)$ counts the number of
times the complement ${\bf 1} -  x$ of $x$ is representable as a linear
combination of $x^1...x^d$. This number is odd ($=1$) iff $x^1...x^d$
are linearly independent and is even otherwise. Thus 
$$
U^d(f_0) = 1 - 2\pr[x^1...x^d\mbox{ are linearly independent}]
$$
and this is easily seen (and well-known) to be a positive constant
bounded away from $1$.

The proof that the acceptance probability of $f$ is high closely
follows the proof of Proposition~15 in \cite{ST00}. That proposition is a
special case $d=2$ of Proposition~\ref{pro:hypertest}. We will repeat
parts of the proof since many definitions need to be modified, and for
completeness, but
will omit proofs of intermediary steps if they are similar.

First, the probability of $f$ to be accepted is 
\begin{equation}
\label{eq:acc_pr}
\frac{1}{2^{|E|}}
\sum_{{\cal S} \subseteq E} \E_{x^1,\ldots,x^k}
\left[ \prod_{T \in {\cal S}} \prod_{i\in T} f(x^i) \cdot
f\left(\sum_{i\in T} x^i\right) \right] 
\end{equation}
In order to simplify this expression, we need to introduce some
notation. 
Let ${\cal F} = \left\{\{1\},...,\{k\}\right\} \cup E :=
\left\{F_1...F_q\right\}$, be 
a family of all the vertices and the edges of $H$, viewed as subsets of
$\{1,\ldots,k\}$. 
Let ${\bf A}$ be a $k\times q$ zero-one matrix whose $q$ columns
are given by $F_1...F_q$, which we view as $0,1$ vectors of length
$k$ (in particular, the first $k$ columns of ${\bf A}$ form the
$k\times k$ identity matrix).  
Let $u_T$, for $T \in E$ be a zero-one
vector of length $q$ 
which is $1$ if $F_i$ is either $T$ or a singleton, corresponding to a
vertex than $T$ passes through; and $0$ otherwise. \\
For ${\cal R} \subseteq 2^{[k]}$, let 
$\E(f,{\cal R}) := \E_{x^1,\ldots,x^k}
\left[ \prod_{R \in {\cal R}} f\left(\sum_{l\in R}
x^l\right) \right].$
\\
Let $U = \mbox{ Span }
\left(u_T\right)_{T\in E}$ be a $t$-dimensional subspace of
${\Z}^q_2$ then, for a boolean $f$, 
(\ref{eq:acc_pr}) is 
$$
\frac{1}{2^{|E|}} \sum_{{\cal S} \subseteq E} \E_{x^1,\ldots,x^k}
\left[ \prod_{i~:~\bigoplus_{T\in {\cal S}} u_T(i) = 1}
f\left(\sum_{l\in F_i} 
x^l\right) \right] = 
$$
$$
\frac{1}{2^t}
\sum_{u = (u(1),\ldots u(q)) \in U} \E_{x^1,\ldots,x^k}
\left[ \prod_{i~:~u(i) = 1} f\left(\sum_{l\in F_i}
x^l\right) \right] =
$$
\begin{equation}
\label{linear}
\frac 1{2^t}
\sum_{u = (u(1),\ldots u(q)) \in U} \E(f,\{ F_i : u(i)=1\}).
\end{equation}

We will show many of the terms $\E(f,\{ F_i : u(i)=1\})$ are $1$.

\begin{definition}
A family ${\cal R} \subseteq 2^{[k]}$ is an ``even cover'', iff every
subset $T \subseteq [k]$, $|T| \le d$ is covered an even number of times
by the sets $R \in {\cal R}$.
\end{definition}

The proofs of the following three lemmas are easily adaptable from the
proofs of the corresponding statements (Lemmas~17,18,19) in \cite{ST00}.

\begin{lemma}
\label{lem.bad_prop}
Let 
${\cal R} \subseteq 2^{[k]}$. If ${\cal R}$ is an even cover, than
$\E (f,{\cal R}) = 1$. For any ${\cal R}$,  $\E (f_n,{\cal R}) \ge
0$. 
\end{lemma}
\begin{proof} Omitted.
\end{proof}

\begin{lemma}
\label{lem.many_bad}
The number of vectors $u\in U$, such that the family ${\cal R} =
\left\{F_i~:~u(i) = 1\right\}$ is an even cover, is at least
$\max\left\{1,2^{t-\sum_{i=2}^d {k \choose i}}\right\}$. 
\end{lemma}

\begin{proof} Omitted.
\end{proof}

\begin{lemma}
\label{lem.dim_ineq}
$$
t \le |E|.
$$
\end{lemma}

\begin{proof} Omitted.
\end{proof}

\noindent Proposition~\ref{pro:hypertest} now follows from
lemma~\ref{lem.bad_prop}, lemma~\ref{lem.many_bad},
\\
lemma~\ref{lem.dim_ineq} and
(\ref{linear}).
\end{proof}

\subsection{Lower Bounds for Arbitrary Test}

We now pass to our more general result, that holds for any test, including
tests with completeness smaller than 1. Recall that a test making $q$-queries
is called {\em non-adaptive} if it makes the $q$ queries simultaneously, instead
of using the answer to the first query to decide how to make the second query, and so on.

\begin{theorem}
\label{th:lower_bnd} let $\cal T$ be a
non-adaptive test that makes $q$ queries and that is a degree $(d-1)$ 
relaxed linearity test with completeness $c$ and soundness $s$
for functions $f:\B^n \to \H$.

Then 
\[ 1-c+s \geq 2^{-q+ \Omega \(\frac{1}{q^{(d-1)/d}}\) } \ . \]
Equivalently, the amortized query complexity of $\cal T$ is 
at least $1 + \Omega\(\frac{1}{q^{(d-1)/d}}\)$.
\end{theorem}
\begin{proof}
We will show that $(1-c) + s$ is at least $\epsilon = \epsilon_H - 
2^{-\Omega(n)}$, where $\epsilon_H$ is the best error achievable by a
hypergraph test after $q$ queries. The theorem will then follow from
Proposition~\ref{pro:hypertest}.

By von Neumann's minimax theorem (also known as Yao's pronciple) it is
enough to construct two distributions $P$ and $Q$ on boolean
functions, $P$ supported on linear functions and $Q$ on functions with
small $U^d$ norm, such that 
such that for any subset $X = \{x^1...x^q\}$ of the boolean
cube the distributions $P'$ and $Q'$ induced by $P$ and $Q$ on $\B^q$
by evaluating a function $f\sim P$ (correspondingly $f\sim Q$) on $X$ are
at most $1 - \epsilon'$ apart. 

For a function $f:\B^n \rightarrow \R$, and an $n\times n$ matrix $A$ over
$G{\mathbb F}_2$ let $f_A$ be given by $f_A(x) = f(Ax)$.

Now, fix a non-zero linear function $\ell$, and 
let $f = f(x_1,\ldots,x_n) := (-1)^{x_1 x_2 \cdots x_d + x_{d+1} \cdots
  x_{2d} + \cdots}$. The distribution $P$ is taken to be uniform over the
functions $\left\{\ell_A~:~A\mbox{ is invertible}\right\}$ and the
distribution $Q$ is taken to be uniform over the 
functions $\left\{f_A~:~A\mbox{ is invertible}\right\}$. 

Observe that $P$ is supported on linear functions. In fact, $P$ is
uniform over all non-zero linear functions. On the other hand $f$ has
a small $U^d$ norm, by Proposition~\ref{pro:hypertest}. The following
lemma shows the same for the functions $f_A$.
\begin{lemma}
For a function $f:\B^n \rightarrow \R$, and an invertible 
$n\times n$ matrix $A$ over $G{\mathbb F}_2$ holds 
$$
U^d(f) = U^d\(f_A\)
$$
\end{lemma}
\begin{proof}
$$
U^d\(f_A\) = \E_{x,x^1,...,x^d} \prod_{S \subseteq [d]} f_A\(x +
\sum_{i\in S} x^i \) = \E_{x,x^1,...,x^d} \prod_{S \subseteq [d]} 
f\(Ax + \sum_{i \in S} Ax^i\) = 
$$
$$
\E_{x,x^1,...,x^d} \prod_{S \subseteq
  [d]} f\(x + \sum_{i\in S} x_i \) = U^d(f)
$$
The third equality follows from the fact that if a $(d+1)$-tuple
$(x,x^1,...,x^d)$ is distributed uniformly in $\B^{n(d+1)}$, then so is
the $(d+1)$-tuple $(Ax,Ax^1,...,Ax^d)$.
\end{proof}

Since the distributions $P$ and $Q$ are invariant under invertible
linear transformations, the induced distributions $P'$ and $Q'$ are
determined by the linear structure of $X$, namely linear dependencies
between $x^1...x^q$. Let the linear rank of the vectors $x^1, \ldots, x^q$
over $G\F_2$ be $k$. We will assume $x^1, \ldots, x^k$ are linearly
independent, and $x^{k+i} = \sum_{j=1}^k a_{ij} x^j$ for $i=1, \ldots, q-k$.

In the remainder of the proof we assume $n \gg k$ and call an event
negligible if its probability is exponentially small in $n$.

The distributions $P'$ and $Q'$ are almost precisely modelled by the
following experiments: Choose $y^1,\ldots, y^k$ independently at random, set
$y^{k+i} = \sum_{j=1}^k a_{ij} y^j$, for $i = 1,\ldots, q-k$, and
return $\ell(y^1),\ldots,\ell(y^q)$ or $f(y^1),\ldots,f(y^q)$
correspondingly. The only caveat comes from the negligible event that $y^1,...,y^k$
are linearly dependent.

For $n \gg k$, the distribution $P'$ is, up to a negligible
probability, given by choosing the first $k$ bits uniformly at random,
and setting the other $q-k$ bits according to specified linear
dependencies. 

Consider the distribution $Q'$. By proposition~\ref{pro:hypertest}
the probability that for $i = 1,\ldots, q-k$ holds $f(y^{k+i}) =
\sum_{j=1}^k a_{ij} f(y^j)$ is at least $\epsilon_H$. Observe that in
this case the $q$-tuple we deal with has no linear contradictions and
thus belongs to the support of $P'$. Call such $q$-tuple
``linear''. We have just proved that, up to a negligible factor,  
$$
\pr\{z~:~z\textrm{ is ``linear''}\} \ge \epsilon_H
$$
We have
$$
\|P' - Q'\| = \frac12 \cdot \sum_z \Big | P'(z) - Q'(z) \Big | = 
\frac12 \cdot \sum_z \left[\max\{P'(z), Q'(z)\} - \min\{P'(z),
    Q'(z)\}\right] \le 
$$
$$
1 - \frac12  \cdot \sum_z \min\{P'(z),
    Q'(z)\} \le 1 - \frac12  \cdot \sum_{\textrm{``linear'' }z} \min\{P'(z),
    Q'(z)\}
$$
To complete the proof of the theorem we observe that, up to negligible
factors, $P'(z) \ge 2^{-k}$ for all ``linear'' $z$, and 
$$
Q'(z) = Q'(z_1,\ldots,z_q) \le \pr\{f(y_1) = z_1,...,f(y_k) = z_k\} =
\prod_{i=1}^k \pr\{f(y_i) = z_i\} = 2^{-k}
$$
Therefore 
$$
\sum_{\textrm{``linear'' }z} \min\{P'(z),Q'(z)\} \ge   
\sum_{\textrm{``linear'' }z} P'(z) - 2^{-\Omega(n)} \ge \epsilon_H -
2^{-\Omega(n)} 
$$
\end{proof}

\section{The PCP Construction}
\label{sec:pcp}

\subsection{The Long-Code Test}

We say that a function $g:\B^n \to \H$ is a {\em codeword of the Long Code} (or, simply,
is a long code) if there exists a coordinate $i$ such that $g(x_1,\ldots,x_n) = (-1)^{x_i}$,
that is, if $g=\chi_{ \{ i \} }$. Note that if $g$ is a long code then there is a coordinate
that has degree-1 influence 1 for $g$. This is the extreme case of large low-degree influence
for a bounded function.

Given a collection of $K$ balanced functions $g_1,\ldots,g_K:\B^n \to \H$, we are interested
in designing a test that distinguishes the following two cases:

\begin{itemize}

\item The functions $g_j$ are all equal to the same long code, that is, for some $i\in [n]$
and for all $j\in [K]$,   $g_j(x)=(-1)^{x_i}$; 
\item The  degree-$d$ cross-influence of the collection $g_1,\ldots,g_K$
is less than  than $\delta$, for some small $\delta$ and large $d$.

\end{itemize}

More formally, we have the following definition.

\begin{definition} [$K$-Function Long Code Test]
A test that is given oracle access to $K$ functions $g_1,\ldots,g_K :\B^n \to \H$
is said to have  soundness $s$ and completeness $c$ if the following conditions hold.

\begin{itemize}

\item If the functions $g_j$  are equal to the same long code, then the test accepts
with probability $\geq c$; 
\item For every $\epsilon>0$ there is a $\tau = \tau(\epsilon)>0$ and
  $d=d(\epsilon)$ such that  
if the test accepts with probability $\geq s+\epsilon$, then there is a
  variable 
of degree-$d$ cross-influence at least $\tau$ for the functions $g_j$.
\end{itemize}
\end{definition}

Let $H = ([t],E)$ be a hypergraph on $t$ vertices.

For $0< \gamma < 1/2$, define the distribution $\mu_\gamma$ over $\B^n$ so that
 $\mu_\gamma (x) = \gamma^{w(x)}(1-\gamma)^{n-w(x)}$, where $w(x)$ is
 the number of ones in $x$. 

The $\gamma$-noisy $H$-test is a $(t+|E|)$-function long code test
defined as follows:

\noindent\fbox{\begin{minipage}{2in}
\begin{program}
$\gamma$-noisy-H-Test $(\{g^a\}_{a\in [t] \cup E} )$\+\\
 choose $x,x_1,\ldots,x_t$ uniformly at random in $\B^n$\\
 for every $i\in [t]$, sample $\eta^{i}$ from $\mu_\gamma$\\
 for every $e\in E$, sample $\eta^e$ from $\mu_\gamma$\\
 \ACCEPT\ if and only if\+\\ $\forall e \in E. \ \prod_{i\in e} g^{(i)}\(\eta^{i} +   x_i\) = g^{(e)}\(\eta^{(e)} + 
\sum_{i \in e} x_i \)$ 
\end{program}\end{minipage}}

\bigskip

\begin{remark} In the definition of $K$-function long-code test, we index functions
by integers $1,\ldots,K$, while in the definition of the Hypergraph Test we let
the given functions be indexed by elements of $[t] \cup E$. We hope  the 
reader is not too bothered by this abuse of notation.
\end{remark}

\subsection{Analysis of the Hypergraph Test}

In this section we prove the following theorem.

\begin{theorem} \label{th:inner}
For every hypergraph $H=([t],E)$, and every $\gamma> 0$, 
the $\gamma$-noisy $H$-test is a $(t+|E|)$-function long code
test with completeness $1- (t+1)\gamma|E|$ and soundness $1/2^{|E|}$.
\end{theorem}

The completeness part is clear.

For the soundness part, as in Section \ref{sec:positive-linear}, we can write the probability
that the test accepts a given set of oracle functions as

\[ \sum_{E'\subseteq E} \frac{1}{2^{|E|}} \E_{x_1,\ldots,x_t, \{ \eta^{i},\eta^{e} \}} 
\prod_{e \in E} \prod_{i\in e} g^{(i)}\(\eta^{i} +x_i\) \cdot
g^{(e)}\(\eta^{(e)} + \sum_{i \in e} x_i \) \]

\[ = \frac {1}{2^{|E|}} + \sum_{E'\subseteq E,~E'\neq \emptyset} \frac{1}{2^|E|} \E_{x_1,\ldots,x_t, \{ \eta^{i},\eta^{e} \}} 
\prod_{e \in E} \prod_{i\in e} g^{(i)}\(\eta^{i} +x_i\) \cdot
g^{(e)}\(\eta^{(e)} + \sum_{i \in e} x_i \) \]

and so if thest accepts with probability at least $2^{-|E|} + \epsilon$ there
is a subset of tests $E'\subseteq E$ such that

\[ \E_{x_1,\ldots,x_t, \{ \eta^{i},\eta^{e} \}} 
\prod_{e \in E'} \prod_{i\in e} g^{(i)}\(\eta^{i} +x_i\) \cdot
g^{(e)}\(\eta^{(e)} + \sum_{i \in e} x_i \) > \epsilon \]

It remains to prove the following lemma ($E$ in the lemma plays the
role of $E'$ above).

\begin{lemma}\label{lm:influence-noisy}
Let $g_j:\B^n \to \H$ be  functions and $H=([t],E)$ be a
hypergraph such that 
\begin{equation}
\E_{x_1,\ldots,x_t, \{ \eta^{i},\eta^{e} \}} 
\prod_{e \in E} \prod_{i\in e} g^{(i)}\(\eta^{i} +x_i\) \cdot
g^{(e)}\(\eta^{(e)} + \sum_{i \in e} x_i \)  >\eps
 \label{eq:test-noisy}
\end{equation}
where the $\eta$'s are sampled according to $\mu_\gamma$.
Then there is a variable $i$  that has degree-$d(\epsilon,\gamma)$
cross-influence at least $\delta (\eps,\gamma)>0$ for the functions $\{ g^{(a)} \}_{a\in [t]\cup E}$.
\end{lemma}

\begin{proof}
Since the $g^{(a)}$ map to $\H$, we can write Equation \ref{eq:test-noisy} equivalently
as

\begin{equation}
 \E_{x_1,\ldots,x_t, \{ \eta^{i},\eta^{e} \}} 
\left(\prod_{i\in {\rm Odd}} g^{(i)}\(\eta^{(i)} +x_i\) \right) \cdot
\left( \prod_{e \in E} g^{(e)}\(\eta^{(e)} \sum_{i \in e} x_i \) \right) > \eps
 \label{eq:test-noisy-equiv}
\end{equation}
Where ${\rm Odd}$ denotes the set of vertices of odd degree in the hypergraph $([t],E)$.

Now, define $G(x) = \E_{\eta} g(\eta + x)$, where $\eta$ is
sampled from $\mu_\gamma$. 

Then 
 (\ref{eq:test-noisy}) becomes 
\begin{equation}
 \E_{x_1,\ldots,x_t} 
\left( \prod_{i\in {\rm Odd}} G^{(i)}\(x_i\) \right) \cdot
\left( \prod_{e \in E} G^{(e)}\(\sum_{i\in e} x_i\) \right) >\eps 
 \label{eq:no-eta}
\end{equation}

Let $k$ be the maximum size of an hyperedge in $E$, and assume,
without loss of generality that $H$ has an hyperedge $e = (1,2,3,\ldots,k)$. Fix the variables   
$x_{k+1},\ldots,x_{t}$ in (\ref{eq:no-eta}) in such a way that the average
over $x_{1},x_{2},\ldots,x_{k}$ is still at least $\epsilon$. In particular,
(\ref{eq:no-eta}) becomes 
\begin{equation} \label{eq:averaged-noisy}
 \E_{x_{1},x_{1},\ldots,x_{k}} f_{\emptyset} (\bzero)
 f_{\{ 1\} } \(x_{1}\)f_{\{ 2 \} }\(x_{2}\) \cdots G^{(e)}\(x_{1} + x_{2} +
 \ldots + x_{k}\) > \eps  
\end{equation}
Where $f_\emptyset (x)$ is the constant function
equal product of all the terms in (\ref{eq:no-eta}) that depend
exclusively on  the fixed variables $x_{(k+1)},\ldots,x_{(t)}$,
$f_{\{ 1\}} \(x_{1}\)$ is the 
product of all the terms that depend only on $x_{1}$ and on
$x_{k+1},\ldots,x_{t}$, and so on. 

Note that each $f_S$ is a product of shifts of functions $G^{(a)}$. Furthermore, the
index sets of $a$'s for distinct $f_S$ are disjoint.

If we call $f_{[k]}:= G^{(e)}$, then Equation \ref{eq:averaged-noisy} says that
\[ \lip{ \{ f_S \} }{k} > \eps \]
Choose $\delta= \delta(k,\eps) =\epsilon^4/2^{O(k)}$ to be small enough so that Lemma \ref{lm:lip} implies
that there are two functions $f_S,f_T$, with $S\neq T$, and a variable $i$, such that
$\I_i (f_S)$ and $\I_i(f_T)$ are both at least $\delta$. 

By previous observations
on the relation between the $f_S$ and the $G^{(a)}$ and by Lemma \ref{lm:influence-shift-general}, we have
that there is a $\delta'$  and two functions $G^a,G^b$
such that $\I_i (G^a)$ and $\I_i (G^b)$ are both at least $\delta'$. 

Consider now the Fourier transform of a function $G(x):= \E_\eta g(\eta+x)$. It
is easy to see that the Fourier coefficients of $G$ satisfy $\hat G(\alpha) = (1-2\gamma)^{|\alpha|}\hat g(\alpha)$.
Therefore, for every degree bound $d$, we have

\[ \I_i (G) = \sum_{\alpha: i\in \alpha} \hat G^2(\alpha) =
\sum_{\alpha: i\in \alpha} (1-2\gamma)^{|\alpha|} \hat g^2(\alpha) \leq
(1-2\gamma)^{d} + \sum_{\alpha: i\in \alpha, |\alpha| \leq d} g^2(\alpha) \leq
(1-2\gamma)^d + \I^{\leq d} (g) \]
This means that we can get $\I^{\leq d}_i\(g^{(a)}\), \I^{\leq d}_i\(g^{(b)}\) \geq
\delta'/2$ if we choose $d= 
O\(\gamma^{-1} \log (\delta'^{-1})\)$.  
\end{proof}

\subsection{Composition and PCP Construction}

The following theorem follows from \autoref{th:outer} and \autoref{th:inner}
using standard techniques.

We need a couple of definitions. If $f:\B^n \to \H$ is a boolean function,
then we define its {\em folding} as the boolean function $\overline f$
defined as follows: $\overline f(0,x_2,\ldots,x_n) := f(0,x_2,\ldots,x_n)$ and
$\overline f(1,x_2,\ldots,x_n) := -f(0,1-x_2,\ldots,1-x_n)$.

The definition satisfies the following useful properties: (i) if $f$ is
a long code, then $f=\overline f$, and, (ii) for every $f$, $\E_x \overline f(x) = 0$, that
is, $\hat { \overline f } (\emptyset) = 0$.

Let $f:\B^n \to \R$ be a  function, and $\pi: [n]\to [n]$ be
a permutation in $S_n$. Then we let $f\circ \pi: \B^n \to \R$ be
the function such that $f\circ \pi (x_1,\ldots,x_n):= f(x_{\pi(1)},\ldots,x_{\pi(n)})$.
Here the interesting properties are that, for every coordinate $i$ and degree bound $t$,
\[ \I^{\leq t}_{\pi(i)} (f) = \I^{\leq t}_{i} (f\circ \pi) \]
and that if $f$ is the long code of $i$, $g$ is the long code of $j$, and $\pi(i) = \pi'(j)$,
then $f\circ \pi$ and $g\circ \pi'$ are the same function.

\begin{theorem}[Main]\label{th:main}
Suppose that the Unique Games Conjecture is true. Then, for every $\delta>0$,
and for every $q\geq 3$,
\[ \np = \pcp_{1-\delta,2q/2^q+\delta} [O(\log n) , q] \]
and, if we can write $q=2^k-1$, then
\[ \np = \pcp_{1-\delta,(q+1)/2^q+\delta} [O(\log n) , q] \]

\end{theorem} 

\begin{proof}
Let $H=([t],E)$ be an hypergraph such that $t+|E| = q$ and $t\leq 1+\log_2 q$. 
(Or $t=\log_2 (q+1)$ if it is possible to write $q=2^t-1$.)

Fix a small constant $\delta$, and consider the $\delta$-noisy $H$-test.
{From} \autoref{th:inner} we have that there are constants
 $\tau(\delta)$ and $d(\delta)$ such that if $\{ g^a \}_{a\in [t] \cup E}$
 are functions accepted by the test with probability at least $1/2^|E|$ then
 the degree-$d(\delta)$ cross-influence of the functions is at least $\tau(\delta)$.
 
Let $\gamma$ be a constant smaller than $\delta \tau^2/d^2$.

{From} \autoref{th:outer} we know that, if the Unique Games conjecture is true,
there is a reduction from SAT to $q$-ary unique games with the property that
a satisfiable instance $\phi$ of SAT is mapped into a unique game $U_\phi$ of strong value 
at least $1-\gamma$ and an unsatisfiable instance $\phi$ of SAT is mapped into a unique
game $U_\phi$ of weak value at most $\gamma$.

A PCP for a formula $\phi$ is a long code for each of the variables of the unique game
$U_\phi$. Let $n$ be size of the alphabet of $U_\phi$.

We consider the following verifier:

\begin{itemize}

\item Pick at random a constraint of $U_\phi$, say that it involves the variables
$v_1,\ldots,v_q$ and the permutations $\pi_1,\ldots,\pi_q$. Let $f_1,\ldots,f_q$
be the functions written in the proof being checked and which, supposedly, are
the long codes of assignments for the variables $v_1,\ldots,v_q$.

\item Run the $\delta$-noisy $H$-test using the functions $\overline f_1\circ \pi_1,\ldots,\overline f_q \circ \pi_q$.

\end{itemize}

\paragraph{Completeness of the test}
If $\phi$ is satisfiable, then consider the proof where for each variable
$v$ we write the function $f$ that is the long code of $A(v)$, with $A()$
being an assignment that strongly satisfies a $1-\gamma$ fraction of constraints.
Then the test accepts with probability at least $(1-\gamma)\cdot (1-q\delta) \geq 1-(q+1)\delta$.

\paragraph {Soundness of the test}
Suppose that the test accepts with probability at least $1/2^{|E|} + 2\delta$.
Consider the following randomized assignment: for every variable $v$, consider
the function $f$ written in the proof corresponding to $v$; consider the
set of coordinates $i$ such that $\I_i ^{\leq d} (\overline f) \geq \tau$; if the
set is empty, give $v$ a random value, otherwise, give $v$ a randomly chosen value
from the set. Note that if the set is non-empty then it has size at most $d/\tau$.

Call a constraint {\em good} if the $H$-test accepts
with probability at least $1/2^{|E|} +\delta$ when that
constraint is chosen by test. Then at least a $\delta$ fraction of constraints is good.

Consider a good constraint, and say that it involves the variables
$v_1,\ldots,v_q$ and the permutations $\pi_1,\ldots,\pi_q$. Let $f_1,\ldots,f_q$
be the functions written in the proof corresponding to the variables $v_1,\ldots,v_q$.

Since the constraint is good, there is a variable $i$ that has degree-$d$ cross-influence 
at least $\tau$ for the 
functions $\overline f_j \circ \pi_j$. That is, there are two functions
$f_a,f_b$ such that $\I^{\leq d}_i (\overline f_a\circ \pi_a) \geq \tau$
and $\I^{\leq d}_i (\overline f_b\circ \pi_b) \geq \tau$. Then the
randomized assignment described above assigns with probability
at least $\tau/d$ the value $\pi_a ^{-1} (i)$ to $v_a$,
and with probability
at least $\tau/d$ the value $\pi_b ^{-1} (i)$ to $v_b$.
When this happens, the constraint is weakly satisfied.

Overall, the randomized assignments weakly satisfies on average
at least a $\delta \tau^2/d^2 > \gamma$ fraction of
constraints, which proves that $\phi$ was satisfiable.

\paragraph{Conclusion} Wrapping up, our PCP verifier has completeness
at least $1- q\delta$ and soundness at most $1/2^|E| + 2\delta$.
Since $\delta$ was chosen arbitrarily at the beginning, the theorem follows.
\end{proof}

\subsection{Inapproximability Results}

It is an immediate consequence of \autoref{th:main} that, assuming the unique games conjecture,
Max $k$CSP cannot be approximated within a factor larger  than $(k+1)/2^k$ if $k$
is of the form $2^t-1$. It follows from
\autoref{th:main} and the reductions in \cite{T01} that, assuming the unique games conjecture,
the Maximum Independent Set problem in graphs of maximum degree $D$ cannot be approximated
within a factor larger $(\log D)^c/D$, for sufficiently large $D$,
where $c$ is an absolute constant.

\section*{Acknowledgements}
We are grateful to Elchanan Mossel for his help with the definitions in Section \ref{sec:complex}
and with the proof of \lemmaref{lm:influence-shift-complex}.

\newcommand{\etalchar}[1]{$^{#1}$}

\end{document}